\def\versionno{ 5ecm   --  v 3.2    -- by jf --   14.1.09   }

\documentclass[12pt]{article}
\usepackage{latexsym,amsmath,amssymb,amsfonts,xypic,bbm}

\usepackage[T1]{fontenc}
\usepackage{amsthm}
\usepackage[all]{xy}
\usepackage{graphicx}

\usepackage{epstopdf,hyperref} 
\usepackage[all]{xy}
\usepackage[mathscr]{eucal}

\numberwithin{equation}{section}

\catcode`\@=11
\newif\if@fewtab\@fewtabtrue
{\count255=\time\divide\count255 by 60
\xdef\hourmin{\number\count255}
\multiply\count255 by-60\advance\count255 by\time
\xdef\hourmin{\hourmin:\ifnum\count255<10 0\fi\the\count255}}
\def\ps@draft{\let\@mkboth\@gobbletwo
    \def\@oddfoot{\hbox to 7 cm{\tiny \versionno
       \hfil}\hskip -7cm\hfil\rm\thepage \hfil {\tiny\draftdate}}
    \def\@oddhead{}
    \def\@evenhead{}\let\@evenfoot\@oddfoot}
\def\draftdate{\number\month/\number\day/\number\year\ \ \ \hourmin }

\catcode`\@=12

\theoremstyle{plain}
\newtheorem{theorem}{Theorem}

\theoremstyle{definition}
\newtheorem{definition}[theorem]{Definition}

\def\complex       {\ensuremath{\mathbbm C}}

\def\R             {\mathbbm{R}}

\def\Z             {\mathbb{Z}}
\def\im            {\mathrm{i}}
\def\id            {\mathrm{id}}

\def\bB            {\bQ}  
\def\be            {\begin{equation}}
\def\bearl         {\begin{array}{l}}
\def\bearll        {\begin{array}{ll}}
\def\bQ            {{\dot{\mathcal B}}}
\def\Bun           {\mathcal{B}\hspace{-.8pt}\mbox{\small\sl un}}
\def\Bunc          {\Bun\!^\nabla\!}
\def\Buntriv       {\Bun\mbox{\small\sl triv}}
\def\Buntrivc      {\Buntriv^{\!\nabla}\!}
\def\C             {\ensuremath{\mathcal C}}
\def\calc          {\ensuremath{\mathcal C}}
\def\cali          {{\cal I}}
\def\caln          {\mathcal{N}}
\def\calg          {{\cal G}}
\def\calt          {{\cal T}}
\def\Cat           {\ensuremath{\mathcal{C}\hspace{-1.1pt}\mbox{\small\sl at}}}

\def\cir           {\,{\circ}\,}
\def\CX            {\ensuremath{C(\X)}}
\def\Cob           {\ensuremath{\mathcal C\hspace{-1.1pt}\mbox{\sl ob}}}
\def\curv          {\mathrm{curv}}
\def\dD            {\dQ}  
\def\Desc          {\mathcal{D}\hspace{-.8pt}\mbox{\small\sl esc}}
\def\dQ            {{\dot{\mathcal D}}}  
\def\dlog          {\mathrm{dlog}\,}
\def\E             {\ensuremath{\mathrm E}}  
\def\ee            {\end{equation}}
\def\eE            {{\rm e}}
\def\eear          {\end{array}}

\def\eq            {\,{=}\,}
\newcommand\erf[1] {(\ref{#1})}
\newcommand\Frac[2]{\mbox{\large$\frac{#1}{#2}$}}
\def\g             {\mathfrak{g}}
\def\GammaM        {\ensuremath{\Gamma_{\!\!\Mf}}}
\def\GammaMX       {\ensuremath{\Gamma_{\!\!\MX}}}

\def\Grbc          {\Grb^{\!\nabla}\!}
\def\Grb           {\mathcal{G}\hspace{-.8pt}\mbox{\small\sl rb}}
\def\Grbctriv      {\Grb\mbox{\small\sl triv}^{\!\nabla}\!}
\def\hol           {\mathrm{Hol}}
\def\Hom           {\mathrm{Hom}}
\def\ii            {{\rm i}}
\def\iN            {\,{\in}\,}
\def\itx           {\item[\raisebox{.08em}{\rule{.44em}{.44em}}]}
\def\Jan           {\mathcal{J}\hspace{-2.9pt}\mbox{\small\sl an}^{\!\nabla}\!}
\def\Jantriv       {\mathcal{J}\hspace{-2.9pt}\mbox{\small\sl antriv}^{\!\nabla}\!}
\def\M             {\ensuremath{\mathcal{M}\hspace{-.8pt}\mbox{\small\sl an}}}
\def\MI            {\M_+}
\def\Mf            {\ensuremath{\mathscr M}}   
\def\Mopp          {\M^{\mathrm{opp}}}
\def\MoppI         {\M^{\mathrm{opp}}_{+}}
\def\Mt            {\ensuremath{M}}  
\def\MX            {\ensuremath{\Mf_\X}}
\def\one           {{\bf1}}
\def\ort           {\ensuremath{\mathrm{or}_2}}
\def\oti           {\,{\otimes}\,}
\def\reals         {{\ensuremath{\mathbb R}}}
\def\rhor          {{\reflectbox{$\rho$}}}

\def\tftC          {\ensuremath{\mbox{\tt tft}_\C}}
\def\Times         {\hspace{.5pt}{\times}\hspace{.5pt}}
\def\To            {\hspace{-.1pt}{\to}\hspace{-.3pt}}
\def\TO            {\hspace{-2.7pt}{\to}\hspace{-2.7pt}}
\def\trivlin       {\mathbf 1}
\def\tr            {\mathrm{tr}}
\def\TTo           {\hspace{-.1pt}{\Rightarrow}\hspace{-.3pt}}
\def\Ue            {\mathrm U(1)}
\def\V             {\ensuremath{\mathscr V}}
\def\vectBc        {\mathcal{V}\hspace{-.8pt}\mbox{\small\sl ect}\Bunc}
\def\X             {\ensuremath{\Sigma}}
\def\Xh            {\ensuremath{\widehat{\Sigma}}}

\makeatletter

\def\adress#1{\gdef\@adress{#1}}
\def\@adress{}
\def\preprint#1{\gdef\@preprint{#1}}
\def\@preprint{}
\def\@maketitle{
  \newpage
  \noindent
  \begin{tabular}{cc}
    \begin{minipage}[c]{0.4\textwidth}
      \begin{flushleft}
        \includegraphics[width=110pt]{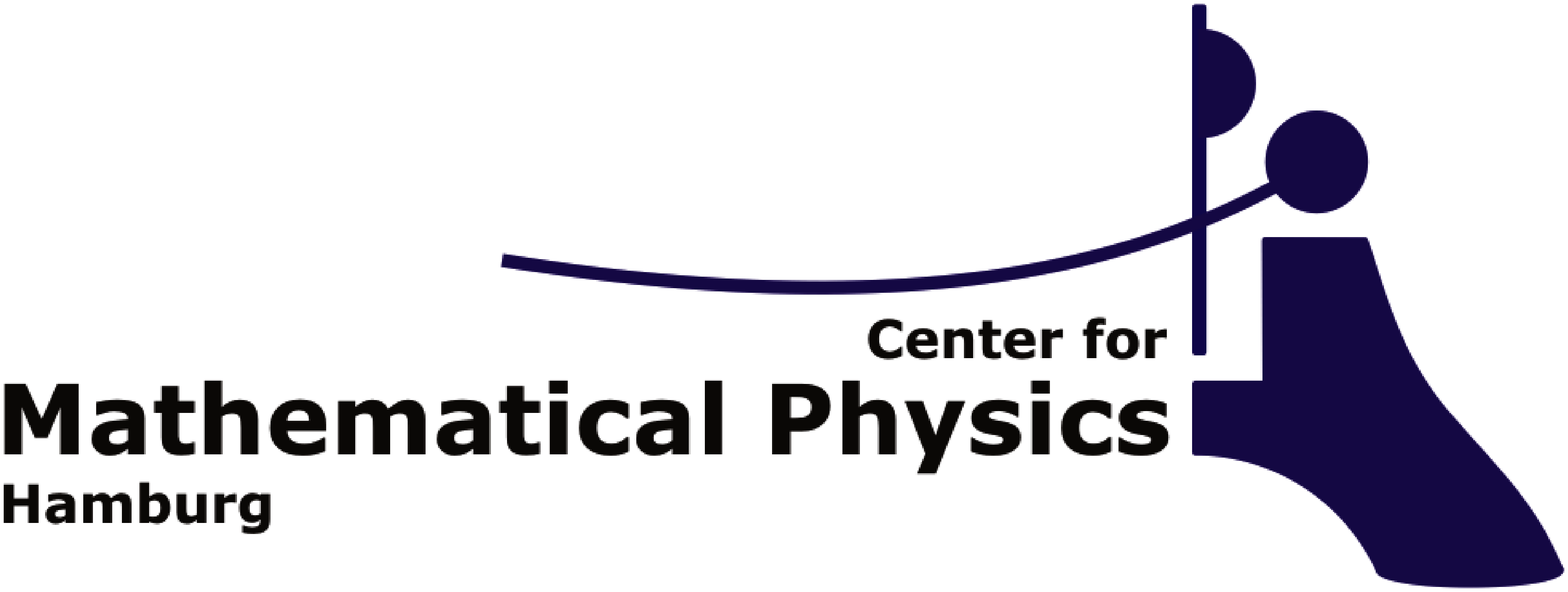}      
      \end{flushleft}
    \end{minipage}&
    \begin{minipage}[c]{0.55\textwidth}
      \begin{flushright}
      {\small\sf\@preprint}
      \end{flushright}
    \end{minipage}
  \end{tabular}
  \vskip 3cm
  \begin{center}
    \LARGE\@title
    \if!\@author!\else \vskip 1cm \large\@author\fi
    \if!\@adress!\else \vskip 1cm \normalsize\@adress\fi
  \end{center}
  \vskip 2cm
}

\newcommand{\alxydim}[2]{\begin{aligned}\xymatrix#1{#2}\end{aligned}}
\newcommand{\bbf}{\varpi}

\makeatother

\def\to{\!\xymatrix@C=0.5cm{\ar[r] &}\!}
\def\iso{\!\xymatrix@C=0.55cm{\ar[r]^-{\sim} &}\!}
\def\mapsto{\!\xymatrix@C=0.5cm{\ar@{|->}[r] &}\!}
\def\Rightarrow{\!\xymatrix@C=0.5cm{\ar@{=>}[r] &}\!}
\def\incl{\!\xymatrix@C=0.5cm{\ar@{^(->}[r] &}\!}
\def\hookrightarrow{\incl}

\def\longrightarrow{\!\xymatrix@C=0.75cm{\ar[r] &}\!}

\def\Schwerpunkt {Bereich}

\begin{document}


\title{%
BUNDLE GERBES AND SURFACE HOLONOMY}
\author{J\"urgen Fuchs$\,^{1}$,~ Thomas Nikolaus$\,^2$,~ 
Christoph Schweigert$\,^2$,~ Konrad Waldorf$\,\,^3$}
\adress{$^{1}$ Teoretisk fysik\\Karlstads Universitet
\\Universitetsgatan 21\\ S\,--\,\,651\,88 Karlstad \\ \bigskip
$^{2}$ Organisationseinheit Mathematik\\ \Schwerpunkt\ 
Algebra und Zahlentheorie\\Universit\"at
Hamburg\\Bundesstra\ss e 55\\D\,--\,20146 Hamburg \\ \bigskip
$^{3}$ Department of Mathematics\\University of California, Berkeley\\
970 Evans Hall \#3840\\Berkeley, CA 94720, USA
}
\preprint{Hamburger Beitr\"age zur Mathematik Nr.\ 323\\ZMP-HH/09-1
\\[-4em]~}
\maketitle

\thispagestyle{empty}

\begin{abstract}
\noindent
Hermitian bundle gerbes with connection are geometric objects for which a 
notion of surface holonomy can be defined for closed oriented surfaces.
We systematically introduce bundle gerbes
by closing the pre-stack of trivial bundle gerbes under descent.
\\
Inspired by structures arising
in a representation theoretic approach to rational conformal field theories, 
we introduce geometric structure that is appropriate to define surface holonomy 
in more general situations: Jandl gerbes for unoriented surfaces, D-bra\-nes
for surfaces with boundaries, and bi-branes for surfaces with defect lines.
\end{abstract}

\thispagestyle{empty}
\newpage
\setcounter{page}{1}

\section{Introduction}

Two-dimensional quantum field theories have been a rich source of relations 
between different mathematical disciplines.  A prominent
class of examples of such theories are the
two-di\-men\-sional rational conformal field theories, which admit a 
mathematically precise description (see \cite{scfr2} for a summary of recent
progress). A large subclass of these also have a classical description
in terms of an action, in which a term given by a surface holonomy enters. 

The appropriate geometric object for the definition of surface holonomies 
for oriented surfaces with empty boundary
are hermitian bundle gerbes. We systematically introduce bundle gerbes by 
first defining a pre-stack of trivial bundle gerbes, in such a way that surface 
holonomy can be defined, and then closing this pre-stack under descent. This
construction constitutes in fact a generalization of the geometry of line 
bundles, their holonomy and their applications to classical particle mechanics. 

Inspired by results in a representation theoretic approach to rational conformal
field theories, we then introduce in the same spirit geometric structure that 
allows to define surface holonomy in more general situations:
Jandl gerbes for unoriented surfaces, D-branes for surfaces with boundaries,
and bi-branes for surfaces with defect lines.


\section{Hermitian line bundles and holonomy}\label{fnsw:sec:bun}

Before discussing bundle gerbes, it is appropriate 
to summarize some pertinent aspects of line bundles.

One of the basic features of a (complex) line bundle $L$ over a smooth manifold 
$M$ is that it is locally trivializable. This means that $M$ can be covered by 
open sets $U_{\alpha}$ such that there exist isomorphisms $\phi_{\alpha}\colon
L|_{U_\alpha} {\to} \trivlin_{U_{\alpha}}$, where $\trivlin_{U_{\alpha}}$ denotes 
the trivial line bundle $ \complex \,{\times}\, U_{\alpha}$. A choice
of such maps $\phi_{\alpha}$ defines gluing isomorphisms
  \begin{equation}
  g_{\alpha\beta}:\quad
  \trivlin_{U_{\alpha}}\big|^{}_{U_\alpha\cap U_\beta}
  \to \trivlin_{U_{\beta}}\big|^{}_{U_\alpha\cap U_\beta} \quad\text{ with }\quad
  g_{\beta\gamma} \circ g_{\alpha\beta} = g_{\alpha\gamma}\,
  ~~\text{on}~~ U_\alpha{\cap} U_\beta{\cap} U_\gamma \,.
  \end{equation}
Isomorphisms between trivial line bundles are just smooth functions. Given a set
of gluing isomorphisms one can 
obtain as additional structure the total space as the manifold
  \begin{equation}
  \label{fnsw:gluing}
  L \,:=\, \bigsqcup_{\alpha} 
  \trivlin_{U_{\alpha}} \,/ \!\sim \mathrm{,}
  \end{equation}
with the relation $\sim$ identifying an element $\ell$ of $\trivlin_{U_{\alpha}}$ 
with $g_{\alpha\beta}(\ell)$ of $\trivlin_{U_{\beta}}$. In short, every bundle 
is glued together from trivial bundles. 

In the following all line bundles will be equipped with a hermitian metric, and 
all isomorphisms are supposed to be isometries. Such line bundles form categories,
denoted $\Bun(M)$. The trivial bundle $\trivlin_M$ defines a full, one-object 
subcategory $\Buntriv(M)$ whose endomorphism set is the monoid of $\Ue$-valued 
functions on $M$. Denoting by $\pi_0(\calc)$ the set of isomorphism classes
of a category $\calc$ and by $H^\bullet(M,\underline{\Ue})$ 
the sheaf cohomology of $M$ with coefficients in the sheaf
of $\Ue$-valued functions, we have the bijection
  \begin{equation}
  \label{fnsw:charclass}
  \pi_0(\Bun(M)) \,\cong\, H^1(M,\underline{\Ue}) \,\cong\, H^2(M,\Z)\,\text{,}
  \end{equation} 
under which the isomorphism class of the trivial bundle is mapped to zero.

Another basic feature of line bundles is that they pull back along smooth maps: for 
$L$ a line bundle over $M$ and $f{:}\ M'\To M$ a smooth map, the pullback $f^*L$ 
is a line bundle over $M'$, and this pullback $f^*$ extends to a functor 
  \begin{equation}
  f^{*}:\quad \Bun(M) \to \Bun(M')\,\text{.}
  \end{equation}
Furthermore, there is a unique isomorphism 
$g^*(f^*L) \To (f \cir g)^*L$ for composable maps $f$ and $g$.

\medskip

As our aim is to discuss holonomies, we should in fact consider a different category,
namely line bundles equipped with (metric) connections.
These form again a category, denoted by $\Bunc(M)$, and this
has again a full subcategory $\Buntrivc(M)$ of trivial 
line bundles with connection. But now this subcategory has more than one object: 
every 1-form $\omega \iN \Omega^1(M)$ can serve as a connection on a trivial line 
bundle $\trivlin$ over $M$;  the so obtained objects are denoted by 
$\trivlin_\omega$. The set $\Hom(\trivlin_\omega,\trivlin_{\omega'})$ of 
connection-preserving isomorphisms $\eta\colon \trivlin_\omega \to \trivlin_{\omega'}$ 
is the set of smooth functions $g\colon M \To \Ue$ satisfying 
  \begin{equation}\label{fnsw:koZyBed}
  \omega' - \omega = -\,{\ii}\, \dlog g\,\text{.}
  \end{equation}
Just like in \erf{fnsw:gluing}, every line bundle $L$ with connection can be glued 
together from line bundles $\trivlin_{\omega_{\alpha}}$ along connection-preserving 
gluing isomorphisms $\eta_{\alpha\beta}$. 

The curvature of a trivial line bundle $\trivlin_\omega$ is 
$\mathrm{curv}(\trivlin_\omega)\,{:=}\,\mathrm{d}\omega\in\Omega^2(M)$, and is 
thus invariant under connection-preserving isomorphisms. It follows that the 
curvature of any line bundle with connection is a globally well-defined, 
closed 2-form. We recall that the cohomology class of this 2-form in real 
cohomology coincides with the characteristic class in \erf{fnsw:charclass}.

In order to introduce the holonomy of  line bundles with connection, we 
say that the holonomy of a trivial line bundle $\trivlin_{\omega}$ over $S^1$ is 
  \begin{equation}
  \hol_{\trivlin_{\omega}}
  := \exp \Big( 2\pi\im \int_{S^1} \omega \Big) \in \Ue\,\text{.}
  \end{equation}
If $\trivlin_{\omega}$ and $\trivlin_{\omega'}$ are trivial line bundles over $S^1$, 
and if there exists a morphism $\eta$ in $\Hom(\trivlin_{\omega},\trivlin_{\omega'})$, 
we have $\hol_{\trivlin_{\omega}} \eq \hol_{\trivlin_{\omega'}}$ because
  \begin{equation}
  \int_{S^1} \omega' - \int_{S^1} \omega
  =  \int_{S^1} -\,\ii\, \dlog\eta ~\in \mathbb{Z}\,\text{.}
  \end{equation}
More generally, if $L$ is any  line bundle with connection over $M$, and 
$\Phi\colon S^1\To M$ is a smooth map, then the pullback bundle 
$\Phi^*L$ is trivial since $H^{2}(S^1,\Z)\eq 0$, and hence one can choose an
isomorphism $\mathcal T\colon \Phi^*L \iso \trivlin_\omega$
for some $\omega\iN\Omega^1(S^1)$. We then set 
  \begin{equation}
  \hol_L(\Phi) := \hol_{\trivlin_{\omega}}\,\text{.}
  \end{equation}
This is well-defined because any other trivialization $\mathcal T'\!\colon 
\Phi^*L \To 1_{\omega'}$ provides a transition isomorphism $\eta \,{:=}\, 
\mathcal T' \cir \mathcal T^{-1}$ in $\Hom(\trivlin_{\omega},\trivlin_{\omega'})$. 
But as we have seen above, the holonomies of isomorphic trivial line bundles coincide.

\medskip

Let us also mention an elementary example of a physical application of line bundles and 
their holonomies: the action functional $S$ for a charged point particle. For 
$(M,g)$ a \mbox{(pseudo-)}Riemannian manifold and 
$\Phi\colon \reals\supset[t_1,t_2] \to (M,g)$ the trajectory of a point particle of
mass $m$ and electric charge $e$, one commonly writes the action $S[\Phi]$ as the sum of
the kinetic term
  \be
  S_{\rm kin}[\Phi] = \frac m2 \int_{t_1}^{t_2}\!
  g( \Frac{{\rm d}\Phi}{{\rm d}t} , \Frac{{\rm d}\Phi}{{\rm d}t} )
  \ee
and a term
  \be
  - e \int_{t_1}^{t_2}\! \Phi^* A \,,
  \ee
with $A$ the electromagnetic gauge potential. However, this formulation is 
inappropriate when the electromagnetic field strength $F$ is not exact, so that a 
gauge potential $A$ with $\mathrm{d}A \eq F$ exists only locally. As explained above, 
keeping track of such local 1-forms $A_{\alpha}$ and  local `gauge transformations', 
i.e.\ connection-preserving isomorphisms between those, leads 
to the notion of a line bundle $L$ with connection. For a closed trajectory, 
i.e.\ $\Phi(t_1) \eq \Phi(t_2)$, the action should be defined as 
  \begin{equation}
  \label{fnsw:action}
  \eE^{\ii S[\Phi]} = \eE^{\ii S_{\rm kin}[\Phi]}\, \hol_L(\Phi)\,\text{.}
  \end{equation}

An important feature of bundles in physical applications is the 
`Dirac quantization' condition on the field strength $F$: the integral 
of $F$ over any closed surface $\Sigma$ in $M$ gives an integer. This follows 
from the coincidence of the cohomology class of $F$ with the characteristic class in
\erf{fnsw:charclass}. Another aspect is a neat explanation of the Aharonov-Bohm effect.
A line bundle over a non-simply connected manifold can have vanishing curvature and yet 
non-trivial holonomies. In the quantum theory holonomies are observable, and thus 
the gauge potential $A$ contains physically relevant 
information even if its field strength is  zero. 
Both aspects, the quantization condition and the Aharonov-Bohm effect, persist in 
the generalization of line bundles to bundle gerbes, which we discuss next.


\section{Gerbes and surface holonomy}

In this section we formalize the procedure of Section \ref{fnsw:sec:bun} that has lead us from 
local 1-form gauge potentials to  line bundles with connection: 
we will explain that it is the closure of 
the category  of trivial bundles with connection under descent. We then apply the same 
principle to locally defined 2-forms, whereby we arrive straightforwardly at the notion 
of bundle gerbes with connection.  We describe the notion of surface holonomy of such 
gerbes and their applications to physics analogously to Section \ref{fnsw:sec:bun}. 


\subsection{Descent of bundles}

As a framework for  structures with a category assigned to every manifold and consistent 
pullback functors we consider presheaves of categories. Let $\M$ be the category 
of smooth manifolds and smooth maps, and let $\Cat$ be the 2-category of 
categories, with functors between categories as 1-morphisms and natural transformations 
between functors as 2-morphisms. Then a presheaf of categories is a lax functor
  \begin{equation}
  \mathcal{F}:\quad \Mopp \to \Cat
  \end{equation}
It assigns to every manifold $M$ a category $\mathcal{F}(M)$, and to every smooth map 
$f\colon M' \To M$ a functor $\mathcal{F}(f)\colon \mathcal{F}(M) \To \mathcal{F}(M')$. 
By the qualification `lax' we mean that the composition of maps must only be preserved
up to coherent isomorphisms. 

In Section \ref{fnsw:sec:bun} we have already encountered four examples of presheaves:
the presheaf $\Bun$ of line bundles, the presheaf $\Bunc$ of line 
bundles with connection, and their sub-presheaves of trivial bundles.

To formulate a gluing condition for presheaves of categories we need to specify 
coverings. Here we choose \emph{surjective submersions} $\pi{:}~ Y\!\To M$. We remark
that every cover of $M$ by open sets $U_{\alpha}$ provides a surjective submersion with 
$Y$ the disjoint union of the $U_{\alpha}$;  thus surjective submersions generalize open 
coverings. This generalization proves to be important for many examples of bundle gerbes, 
such as the lifting of bundle gerbes and the canonical bundle gerbes of compact simple 
Lie groups.

With hindsight, a choice of coverings endows the category $\M$ with a Grothendieck 
topology. Both surjective submersions and open covers define a Grothendieck topology,
and since every surjective submersion 
allows for local sections, the resulting two Grothendieck topologies are equivalent.
And in fact the submersion topology is the maximal one equivalent to open coverings.

Along with a covering $\pi{:}~ Y\!\To M$ there comes a simplicial manifold
  \begin{equation}
  \xymatrix{
  {\cdots\,}  
  \ar@<1.3ex>[r]^{\partial_0}
  \ar@<0.3ex>[r]
  \ar@<-0.7ex>[r]
  \ar@<-1.7ex>[r]_{\partial_3}
  &
  {\,Y^{[3]}\,}
  \ar@<0.9ex>[r]^{\partial_0}
  \ar@<-0.1ex>[r]
  \ar@<-1.1ex>[r]_{\partial_2}
  &
  {\,Y^{[2]}\,}
  \ar@<0.3ex>[r]^{\partial_0}
  \ar@<-0.7ex>[r]_{\partial_1} 
  &
  {\,Y} 
  \ar[r]^{\pi}
  &
  M\,.
  }
  \end{equation}
Here $Y^{[n]}$ denotes the $n$-fold fibre product of $Y$ over $M$,
  \begin{equation}
  Y^{[n]} := \{(y_0,\ldots,y_{n-1}) \iN Y^n \,|\, \pi(y_0)\,{=}\ldots{=}\,\pi(y_{n-1}) \} \,,
  \end{equation}
and the map $\partial_i\colon Y^{[n]} \To Y^{[n-1]}$ omits the $i$th entry. In particular
$\partial_0\colon Y^{[2]} \To Y$ is the projection to the second factor and 
$\partial_1\colon Y^{[2]} \To Y$ the one to the first. All fibre products $Y^{[k]}$ are 
smooth manifolds, and all maps $\partial_i$ are smooth. Now let $L$ be a line bundle 
over $M$. By pullback along $\pi$ we obtain:
\def\leftmargini{3.21em}\\[-1.5em]
\begin{itemize}
\item[(BO1)] An object $\tilde L := \pi^*L$ in $\Bun(Y)$.  
\item[(BO2)] A morphism 
  \begin{equation}
  \phi:\quad \partial_0^*\tilde L\cong \partial_0^* \pi^* L
  \iso
  \partial_1^* \pi^* L \cong \partial_1^*\tilde L
  \end{equation}
in $\Bun(Y^{[2]})$ induced from the identity $\pi\cir \partial_0\eq \pi\cir \partial_1$.
\item[(BO3)] 
A commutative diagram
  \begin{equation}
  \alxydim{@C=0.28cm}{
  \partial_1^*\partial_0^* \tilde L 
\ar@/_2pc/_{\partial_1^*\phi}[rrrrrrr]
  \ar@{=}[r] 
  & 
  \partial_0^*\partial_0^* \tilde L\,
  \ar^-{\partial_0^*\phi}[rr]
  &&
  \,\partial_0^*\partial_1^* \tilde L \ar@{=}[r] &
  \partial_2^*\partial_0^* \tilde L\, \ar^-{\partial_2^*\phi}[rr] &&
  \,\partial_2^*\partial_1^* \tilde L 
  \ar@{=}[r]
  &
  \partial_1^*\partial_1^* \tilde L
  }
  \end{equation}
of morphisms in $\Bun(Y^{[3]})$; or in short, an equality
$\partial_2^* \phi\cir \partial_0^* \phi \eq \partial_1^* \phi$.
\end{itemize}

We call a pair $(\tilde L, \phi)$ as in (BO1) and (BO2) which satisfies
(BO3) a \emph{descent object} in the presheaf $\Bun$. Analogously we obtain for a 
morphism $f\colon L \To L'$ of line bundles over $M$
\def\leftmargini{3.38em}\\[-1.5em]
\begin{itemize} 
\item[(BM1)] 
A morphism $\tilde f := \pi^*f \colon \tilde L \To \tilde L'$ in $\Bun(Y)$.
\item[(BM2)] A commutative diagram
  \begin{equation}
  \phi' \circ \partial_0^*\tilde f = \partial_1^*\tilde f \circ  \phi
  \end{equation}
of morphisms in $\Bun(Y^{[2]})$. 
\end{itemize}
Such a morphism $\tilde f$ as in (BM1) obeying (BM2) is called a \emph{descent morphism} in 
the presheaf $\Bun$. 

Descent objects and descent morphisms for a given covering $\pi$ form a category
$\Desc(\pi{:}\,Y\TO M)$ of descent data. What we described above is a functor
  \begin{equation}
  \iota_\pi:\ \Bun(M) \to \Desc(\pi{:}\,Y\TO M)\,\text{.}
  \end{equation}
The question arises whether every `local' descent object corresponds to a `global' object on $M$,
i.e.\ whether the functor $\iota_{\pi}$ is an equivalence of categories. 

The construction generalizes straightforwardly to any presheaf of categories $\mathcal{F}$, 
and if the functor $\iota_{\pi}$ is an equivalence for all coverings $\pi\colon Y \!\To M$, 
the presheaf $\mathcal{F}$ is called  a \emph{sheaf of categories} (or stack). Extending 
the gluing process from \erf{fnsw:gluing} to non-trivial bundles shows that
the presheaves $\Bun$ and $\Bunc$ are  sheaves. In contrast, the presheaves 
$\Buntriv$ and $\Buntrivc$ of trivial bundles are not sheaves, since
gluing of trivial bundles does in general not result in a trivial bundle. 
In fact the gluing process \erf{fnsw:gluing} shows that every bundle can be obtained by 
gluing trivial ones. In short, the sheaf $\Bunc$ of line bundles with connection is
obtained by closing the presheaf $\Buntrivc$ under descent.


\subsection{Bundle gerbes}
\label{fnsw:sec:defgrb}

Our construction of line bundles started from trivial line bundles
with connection which are just 1-forms on $M$, and 
the fact that 1-forms can be integrated along curves has lead us to the notion of 
holonomy. To arrive at a notion of \emph{surface} holonomy, we now consider a 
category of 2-forms, or rather a 2-category:
\def\leftmargini{1.12em}~\\[-1.5em]
\begin{itemize}\addtolength\itemsep{-5pt}
\itx 
An object is a 2-form $\omega\iN\Omega^2(M)$, called a
\emph{trivial bundle gerbe with connection} and denoted by $\mathcal{I}_{\omega}$.
\itx 
A 1-morphism $\eta\colon \omega \To \omega'$
is a 1-form $\eta\in\Omega^1(M)$ such that $\mathrm{d}\eta=\omega'-\omega$.
\itx 
A 2-morphism $\phi\colon \eta \TTo \eta'$
is a smooth function $\phi: M \to \Ue$ such that $-\ii \,\mathrm{dlog}(\phi)=\eta'-\eta$.
\end{itemize}

There is also a natural pullback operation along maps, induced by pullback on differential forms. 
The given data can be rewritten as a presheaf of 2-categories, as there is a 2-category 
attached to each manifold. This presheaf should now be closed under descent to obtain 
a sheaf of 2-categories. As a first step we complete the morphism categories  
under descent. Since these are categories of trivial line bundles with connections, we set
  \begin{equation}
  \mathrm{Hom}(\cali_\omega,\cali_{\omega'}) 
  := \Bunc_{\!\omega'-\omega}(M)\,\text{,}
  \end{equation}
the category of hermitian line bundles with connection of fixed curvature $\omega'-\omega$. 
The horizontal composition is given by the tensor product in the category of bundles.
Finally, completing the 2-category under descent, we get the definition of a bundle gerbe: 
\begin{definition}\label{fnsw:def:grb}
A bundle gerbe $\calg$ (with connection) over $M$ consists of the following data:
a covering $\pi\colon Y \!\To M$, and for the associated simplicial manifold
  \begin{equation}
  \xymatrix{\,Y^{[4]\,} \ar[r]<.6em>\ar[r]<.2em>\ar[r]<-.2em>\ar[r]<-.6em> &
  \,Y^{[3]\,} \ar[r]<.4em>\ar[r]<.0em>\ar[r]<-.4em>
  &\,Y^{[2]\,} \ar_-{\partial_1}[r]<-.2em> \ar^-{\partial_0}[r]<.2em> 
  &\,Y \ar^-\pi[r] &\,M 
  }
  \end{equation}
\def\leftmargini{3.31em}
\begin{itemize}
\item[(GO1)]
an object $\mathcal{I}_{\omega}$ of $\Grbctriv(Y)$: a 2-form $\omega\in\Omega^2(Y)$;
\item[(GO2)] 
a 1-morphism 
  \begin{equation}
  L\colon~ \partial_0^{*}\mathcal{I}_{\omega} \To \partial_1^{*}\mathcal{I}_{\omega}
  \end{equation} 
in $\Grbctriv(Y^{[2]})$: a line bundle $L$ with connection over $Y^{[2]}$;
\item[(GO3)]
a 2-isomorphism 
  \begin{equation}
  \mu \colon~ \partial_{2}^*L \oti \partial_{0}^*L \,{\Rightarrow}\, \partial_{1}^*L
  \end{equation} 
in $\Grbctriv(Y^{[3]})$: a connection-preserving morphism of line bundles over $Y^{[3]}$;
\item[(GO4)] 
an equality
  \begin{equation}
  \partial_2^*\mu \circ (\id \oti \partial_0^* \mu) = \partial_1^*\mu \circ (\partial_3^*\mu \oti \id)\,
  \end{equation}
of 2-morphisms in $\Grbctriv(Y^{[4]})$.
\end{itemize}
\end{definition}

For later applications it will be necessary to close the morphism categories under a 
second operation, namely direct sums. Closing the category of line bundles with connection 
under direct sums leads to the category of complex vector bundles with connection, i.e.\ we set
  \be
  \label{defhomvect}
  \mathrm{Hom}(\cali_\omega,\cali_{\omega'}) := \vectBc_{\!\omega'-\omega}(M) \,,
  \ee
where the curvature of these vector bundles is constrained to satisfy
  \be
  \frac 1n\, \tr(\curv(L)) = \omega'-\omega \,,
  \ee
with $n$ the rank of the vector bundle. Notice that this does not affect the 
definition of a bundle gerbe, since the existence of the 2-isomorphism $\mu$ 
restricts the rank of $L$ to be one. 

As a next step, we need to introduce 1-morphisms and 2-morphisms between bundle 
gerbes. 1-morphisms have to compare two bundle gerbes $\calg$ and $\calg'$. We 
assume first that both bundle gerbes have the same covering $Y \!\To M$. 

\begin{definition} 
i) A 1-morphism between bundle gerbes $\calg \eq (Y,\omega,L,\mu)$ and 
$\calg' \eq (Y,\omega',L',\mu')$ over $M$ with the same surjective submersion 
$Y \!\To M$ consists of the following data on the associated simplicial manifold
  \be
  \xymatrix{\,Y^{[4]\,} \ar[r]<.6em>\ar[r]<.2em>\ar[r]<-.2em>\ar[r]<-.6em> &
  \,Y^{[3]\,} \ar[r]<.4em>\ar[r]<.0em>\ar[r]<-.4em> &
  \,Y^{[2]\,} \ar_-{\partial_1}[r]<-.2em> \ar^-{\partial_0}[r]<.2em> &
  \,Y \ar^-\pi[r] &\, M \,.
  }
  \ee
\def\leftmargini{3.91em}
\begin{itemize}
\item[(G1M1)] 
a 1-morphism $A\colon \mathcal{I}_{\omega} \to \mathcal{I}_{\omega'}$ in $\Grbctriv(Y)$:
a rank-$n$ hermitian vector bundle  $A$ with connection of curvature 
$\frac 1 n\, \tr(\mathrm{curv}(L)) \eq \omega' - \omega$;
\item[(G1M2)]
a 2-isomorphism $\alpha\colon L' \oti \partial_0^*A {\Rightarrow} \partial_1^*A \oti L$ 
in $\Grbctriv(Y^{[2]})$: a connection-preserving morphism of hermitian vector bundles; 
\item[(G1M3)]
a commutative diagram
  \begin{equation}
  (\id \oti \mu') \circ (\partial_2^*\alpha \oti \id) \circ (\id \oti \partial_0^*\alpha)
  = \partial_1^* \alpha \circ (\mu \oti id)
  \end{equation}
of 2-morphisms in $\Grbctriv(Y^{[3]})$.
\end{itemize}
ii) A 2-morphism between two such 1-morphisms $(A,\alpha)$ and $(A',\alpha')$ consists of
\def\leftmargini{3.85em}\\[-1.5em]
\begin{itemize}
\item[(G2M1)]
a 2-morphism $\beta\colon A {\Rightarrow} A'$ in $\Grbctriv(Y)$:
a connection-preserving morphism of vector bundles;
\item[(G2M2)]
a commutative diagram
  \begin{equation}
  \alpha' \circ (\id \oti \partial_0^*\beta) = (\partial_1^*\beta \oti \id) \circ \alpha
  \end{equation}
of 2-morphisms in $\Grbctriv(Y^{[2]})$.
\end{itemize}
\end{definition}

Since 1-morphisms are composed by taking tensor products of 
vector bundles, a 1-morphism is invertible if and only if 
its vector bundle is of rank one. 

\medskip

In order to define 1-morphisms and 2-morphisms between bundle gerbes with possibly 
different coverings $\pi\colon Y \!\To M$ and $\pi'\colon Y' \!\To M$, we
pull all the data back to a common refinement of these coverings and compare them 
there. We call a covering $\zeta\colon Z \To M$ a common refinement of $\pi$ and $\pi'$ 
iff there exist  maps $s\colon Z \To Y$ and $s'\colon Z \To Y'$ such that
  \be
  \xymatrix{
  Y \ar[dr]_-{\pi}  & Z \ar[l]_-{s} \ar[r]^{s'}\ar[d]_{\zeta}& Y'\ar[dl]^-{\!\!\pi'} \\
  &M&
  }
  \ee
commutes. An example of such a common refinement is the fibre product  
$Z \,{:=}\, Y {\times_M}\, Y' \To M$, with the maps $Z \To Y$ and 
$Z \To Y'$ given by the projections. The important point about a common 
refinement $Z \To M$  is that the maps $s$ and $s'$ induce simplicial maps 
  \begin{equation}
  \alxydim{}{Y^{\bullet} & Z^{\bullet} \ar[l] \ar[r] & Y^{\prime\bullet}}\!\text{.}
  \end{equation}
For bundle gerbes $\calg \eq (Y,\omega,L,\mu)$ and $\mathcal{G}'=(Y',\omega',L',\mu')$ 
we obtain new bundle gerbes with surjective submersion $Z$ by pulling back all the 
data along the simplicial maps $s$ and $s'$. Explicitly,  
$\calg_Z \,{:=}\, (Z,s_0^*\omega, s_1^*L, s_2^* \mu)$ and 
$\mathcal{G}'_Z=(Z,s_0^{\prime *}\omega',s_1^{\prime *}L',s_2^{\prime *}\mu')$. 
Also morphisms can be refined by pulling them back.

\begin{definition}
i) A \emph{1-morphism} between bundle gerbes $\calg \eq (Y,\omega,L,\mu)$ and 
$\calg' \eq (Y',\omega',L',\mu')$ consists of a common refinement $Z \To M$ of the 
coverings $Y \!\To M$ and $Y'\!\To M$ 
and a morphism $(A,\alpha)$ of the two refined gerbes $\calg_Z$ and $\calg'_Z$. 
\\[2pt]
ii) A 2-morphism between 1-morphisms $\mathfrak{m} \eq (Z,A,\alpha)$ and 
$\mathfrak{m}' \eq (Z',A',\alpha')$ consists of a common refinement $W\!\To M$ of 
the coverings $Z \To M$ and $Z' \To M$ (respecting the projections to $Y$ and $Y'$, 
respectively) and a 2-morphism $\beta$ of the refined morphisms 
$\mathfrak{m}_W$ and $\mathfrak{m}_W'$. In addition two such 2-morphisms
$(W,\beta)$ and $(W',\beta')$ must be identified iff there exists a further 
common refinement $V \!\To M$ of $W \!\To M$ and $W' \To M$, compatible with 
the other projections, such that the refined 2-morphisms agree on $V$. 
\end{definition}

For a gerbe $\calg \eq (Y,\omega,L,\mu)$ and a refinement $Z \To M$ of $Y$ the 
refined gerbe $\calg_Z$ is isomorphic to $\calg$. This implies that every gerbe is 
isomorphic to a gerbe defined over an open covering $Z \,{:=}\, \bigsqcup_{i \in I} U_i$. 
Furthermore we can choose the covering in such a way that the line bundle over double 
intersections is trivial as well. When doing so we obtain the familiar description 
of gerbes in terms of local data, reproducing formulas by \cite{alvaO5,gawe1}.
Extending this description to morphisms it is straightforward 
to show that gerbes are classified by the so-called  Deligne cohomology $H^k(M,\mathcal{D}(2))$ in degree two:
  \begin{equation}
  \pi_0(\Grbc(M)) \,\cong\, H^2(M,\mathcal{D}(2)) \,.
  \end{equation} 
Analogously we get the classification of gerbes without connection as 
  \begin{equation}
  \pi_0(\Grb(M)) \,\cong\, H^2(M,\underline{\Ue}) \,\cong\, H^3(M,\Z)\,\text{.}
  \end{equation}


\subsection{Surface holonomy}\label{fnsw:subsec:surfhol}

The holonomy of a trivial bundle gerbe $\mathcal{I}_{\omega}$ over a 
closed oriented  surface $\Sigma$ is by definition
  \begin{equation}
  \hol_{\mathcal{I}_{\omega}} := \exp \!\Big( 2\pi\im\int_{\Sigma} \omega \Big)
  \,\in \Ue\text\,{.}
  \end{equation}
If $\mathcal{I}_{\omega}$ and $\mathcal{I}_{\omega'}$ are two trivial bundle gerbes over 
$\Sigma$ such that there exists a 1-isomorphism $\mathcal{I}_{\omega} \to I_{\omega'}$, 
i.e. a vector bundle $L$ of rank one, we have an equality 
$\hol_{\mathcal{I}_{\omega}} \eq \hol_{\mathcal{I}_{\omega'}}$ because
  \begin{equation}
  \int_{\Sigma} \omega' - \int_{\Sigma}\omega = \int_{\Sigma}\mathrm{curv}(L)
  \,\in \Z\,\text{.}
  \end{equation}

More generally, consider a bundle gerbe $\mathcal{G}$ with connection over 
a smooth manifold $M$, and a smooth map 
  \begin{equation}
  \Phi:\quad \Sigma \to M
  \end{equation}
defined on a closed oriented surface $\Sigma$. Since $H^3(\Sigma,\Z)=0$, the pullback 
$\Phi^{*}\mathcal{G}$ is isomorphic to a trivial bundle gerbe. Hence one can choose a 
trivialization, i.e.\ a 1-isomorphism
  \begin{equation}
  \calt:\quad \Phi^*\mathcal{G} \iso \mathcal{I}_{\omega}
  \end{equation}
and define the holonomy of $\mathcal{G}$ around $\Phi$ by
  \begin{equation}
  \hol_{\calg}(\Phi) := \hol_{\mathcal{I}_{\omega}}\,\text{.}
  \end{equation}
In the same way as for the holonomy of a line bundle with connection, this definition 
is independent of the choice of the 1-isomorphism $\mathcal{T}$. Namely, if 
$\mathcal{T}'\!\colon \Phi^{*}\mathcal{G} \iso \mathcal{I}_{\omega'}$ is 
another trivialization, we have a transition isomorphism 
  \begin{equation}
  \label{fnsw:transiso}
  L:=\calt'\circ\calt^{-1}\colon~~\cali_\omega \iso \cali_{\omega'} \,,
  \end{equation}
which shows the independence.


\subsection{Wess-Zumino terms} \label{fnsw:sec-wz}

As we have seen in Section \ref{fnsw:sec:bun}, the holonomy of a line bundle with 
connection supplies a term in the action functional of a classical charged particle, 
describing the coupling to a gauge field whose field strength is the curvature of 
the line bundle. Analogously, the surface holonomy of a bundle gerbe with 
connection defines a term in the action of a classical charged string. Such a 
string is described in terms of a smooth map $\Phi\colon \Sigma \To M$.
The exponentiated action functional of the string is (compare \erf{fnsw:action})
  \begin{equation}
  \eE^{\ii S[\Phi]} = \eE^{\ii S_{\rm kin}[\Phi]}\, \hol_\mathcal{G}(\Phi)\,,
  \end{equation}
where $S_{\mathrm{kin}}[\Phi]$ is a kinetic term which involves a conformal structure 
on $\Sigma$. Physical models whose fields are maps defined on surfaces are called 
(non-linear) \emph{sigma models}, and the holonomy term is called a 
\emph{Wess-Zumino term}.  Such terms are needed in certain models in order 
to obtain quantum field theories that are conformally invariant.
 
A particular class of sigma models with Wess-Zumino term is given by WZW
(Wess-Zumino-Witten) models. For these the target space $M$ is a connected 
compact simple Lie group $G$, and the curvature of the bundle gerbe 
$\mathcal{G}$ is an integral  multiple of the canonical 3-form
  \begin{equation*}
  H=\left \langle  \theta \wedge [\theta \,{\wedge}\, \theta] \right \rangle
  \,\in \Omega^3(G)
  \end{equation*}  
($\theta$ is the left-invariant Maurer-Cartan form on $G$, and
$\langle\cdot\,,\cdot\rangle$ the Killing form of the Lie algebra $\g$ of $G$). 
WZW models have been a distinguished arena for the interplay between Lie theory
and the theory of bundle gerbes \cite{gawe1,gawedzki1}. This has lead to
new insights both in the physical applications and in the underlying mathematical 
structures. Some of these will be discussed in the following sections. 

\medskip

Defining Wess-Zumino terms as the holonomy of a bundle gerbe with connection
allows one in particular to explain the following two facts. 
\begin{itemize}
\itx
{\em The Aharonov-Bohm effect\/}: This occurs when the bundle gerbe has a \emph{flat} connection, 
i.e.\ its curvature $H \iN \Omega^3(M)$ vanishes. This does not mean, though, that the 
bundle gerbe is trivial, since its class in $H^3(M,\Z)$ may be pure torsion. In 
particular, it can still have non-constant holonomy, and thus a non-trivial 
Wess-Zumino term. 
\\[2pt]
An example for the Aharonov-Bohm effect is the sigma model on the 2-torus 
$T \eq S^1 \,{\times}\, S^1$. By dimensional reasons, the 3-form $H$ vanishes. Nonetheless, 
since $H^2(T,\Ue)\eq\Ue$, there exists a whole family of Wess-Zumino terms parameterized 
by an angle, of which only the one with angle zero is trivial.

\itx
{\em Discrete torsion\/}: The set of isomorphism classes of bundle gerbes with connection 
that have the same curvature $H$ is parameterized by $H^2(M,\Ue)$ via the map
  \begin{equation}
  H^2(M,\Ue) \to \mathrm{Tors}(H^3(M,\Z))\,\text{.}
  \end{equation}
If this group is non-trivial, there exist \emph{different} Wess-Zumino terms for 
one and the same field strength $H$; their difference is called `discrete torsion'. 
\\[2pt]
An example for discrete torsion is the level-$k$ WZW model on the 
Lie group $\mathrm{PSO}(4n)$. Since $H^2(\mathrm{PSO}(4n),\Ue)\eq\Z_2$, there 
exist two non-isomorphic bundle gerbes with connection having equal curvature. 
\end{itemize}


\section{The representation theoretic formulation of RCFT} \label{fnsw:sec4}

\subsection{Sigma models}

Closely related to surface holonomies are novel geometric structures that
have been introduced for unoriented surfaces, for surfaces with boundary,
and for surfaces with defect lines. These structures constitute the second 
theme of this contribution, extending the construction of gerbes and 
surface holonomy via descent; they will be discussed
in Sections \ref{fnsw:sec5}, \ref{fnsw:sec6} and \ref{fnsw:sec7}.
   
These geometric developments were in fact strongly inspired by algebraic and
representation theoretic results in two-dimensional 
quantum field theories. To appreciate this connection we briefly review 
in this section the relation between spaces of maps 
$\Phi\colon \X\To\Mt$, 
as they appear in the treatment of holonomies, and quantum field theories.

As already indicated in Section \ref{fnsw:sec-wz}, a classical field theory, the
(non-linear) \emph{sigma model}, on a two-dimensional surface \X, called the 
\emph{world sheet}, can be associated to the space of smooth maps 
$\Phi$ from \X\ to some smooth manifold \Mt, called the \emph{target space}.
Appropriate structure on the target space determines a Lagrangian for the
field theory on \X. Geometric structure on \Mt, e.g.\ a (pseu\-do-)\,Rie\-mannian 
metric $G$, becomes, from this point of view, for any given map $\Phi$ a background
function $G(\Phi(x))$ for the field theory on \X.

Three main issues will then lead us to a richer structure related to
surface holonomies:
\def\leftmargini{1.12em}~\\[-1.5em]
\begin{itemize}
\itx 
In string theory (where the world sheet \X\ arises as the surface swept out by 
a string moving in \Mt) and in other applications as well, one also encounters 
sigma models on world sheets \X\ that have {\em non-empty boundary\/}.
We will explain how the geometric data relevant for encoding boundary
conditions -- so called D-branes -- can be derived from geometric principles.
\itx 
String theories of type I, which form an integral part of string dualities, 
involve {\em unoriented\/} world sheets. In string theory it is therefore
a fundamental problem to exhibit geometric structure
on the target space that provides a notion of holonomy for unoriented surfaces.
\itx 
An equally natural structure present in quantum field theory are 
{\em topological defect lines\/},
along which correlation functions of bulk fields can have a branch-cut.
In specific models these can be understood, just like boundary conditions, as
continuum versions of corresponding structures in lattice models of
statistical mechanics. (For instance, in the lattice version of the
Ising model a topological defect is produced by changing the coupling along 
all bonds that cross a specified line from ferromagnetic to antiferromagnetic.)
\end{itemize}

Sigma models have indeed been a significant source of examples for quantum field 
theories, at least on a heuristic level. Conversely, 
having a sigma model interpretation for a given quantum field theory allows 
for a geometric interpretation of quantum field theoretic quantities.
Sigma models have indeed been, at least on a heuristic level, a significant 
source of examples for quantum field theories. Conversely, quantum
field theoretic structures in sigma models have lead to
structural insights and quantitative predictions in geometry.

A distinguished subclass of theories in which this relationship between
quantum field theory and geometry can be studied are two-dimensional conformal 
field theories, or CFTs, for short, and among these in particular the
rational conformal field theories for which there exists a rigorous
representation theoretic approach. The structures appearing in that approach in 
the three situations mentioned above suggest new geometric notions for conformal 
sigma models. Below we will investigate these notions with the help of standard 
geometric principles. Before doing so we formulate, in representation theoretic
terms, the relevant aspects of the quantum field theories in question.


\subsection{Rational conformal field theory}

The conformal symmetry, together with further, so-called chiral, symmetries 
of a CFT can be encoded in the structure of a conformal vertex algebra \V. 
For any conformal vertex algebra one can construct (see e.g.\ \cite{FRbe2}) 
a chiral CFT; in mathematical terms, a chiral CFT is a system of 
{\em conformal blocks\/}, i.e.\ sheaves over the moduli spaces of curves with 
marked points. These sheaves of conformal blocks are endowed with a projectively 
flat connection, the Knizhnik-Za\-mo\-lod\-chikov connection, which in turn
furnishes representations of the fundamental groups of the moduli spaces,
i.e.\ of the mapping class groups. 

Despite the physical origin of its name, a chiral conformal field theory is 
mathematically rigorous. On the other hand, from the two-dimensional point of 
view it is, despite its name, not a conventional quantum field theory, as 
one deals with (sections of) bundles instead of local correlation functions.
In particular, it must not be confused with a full local conformal field
theory, which is the relevant structure to enter our discussion of holonomies.

Chiral conformal field theories are particularly tractable when the 
vertex algebra \V\ is rational in the sense of \cite[thm\,2.1]{huan21}. 
Then the representation category \C\ of \V\ is a modular tensor category,
and the associated chiral CFT is a {\em rational chiral CFT\/}, or chiral RCFT. 
In this situation, we can use the tools of three-dimensional topological 
quantum field theory (TFT). A TFT is, in short, a monoidal functor \tftC\ 
\cite[chap.\,IV.7]{TUra} that associates a finite-dimensional
vector space $\tftC(\E)$ to any (extended) surface \E, and a linear map
from $\tftC(\E)$ to $\tftC(\E')$ to any (extended) cobordism $\Mf\colon \E\To\E'$.

More precisely, a three-dimensional TFT is a projective monoidal functor from
a category $\Cob_\C$  of decorated cobordisms to the category of 
finite-dimensional complex vector spaces. The modular tensor category \C\ 
provides the decoration data for $\Cob_\C$. Specifically, the objects \E\ 
of $\Cob_\C$ are extended surfaces, i.e.\,%
  \footnote{~Here various details are suppressed. 
  Detailed information, e.g.\ the precise definition of a ribbon graph or
  the reason why \tftC\ is only projective, can be found in many places,
  such as \cite{TUra,BAki,KArt} or \cite[sect.\,2.5-2.7]{fffs3}.} 
compact closed oriented two-manifolds with a
finite set of embedded arcs, and each of these arcs is marked by an
object of \C. A morphism $\E \To \E'$ is an {\em extended cobordism\/},
i.e.\ a compact oriented three-manifold \Mf\ with
$\partial\Mf\eq(-\E)\,{\sqcup}\,\E'$, together with an oriented ribbon graph
\GammaM\ in \Mf\ such that at each marked arc of $(-\E)\,{\sqcup}\,\E'$
a ribbon of \GammaM\ is ending. Each ribbon of \GammaM\ is labeled by an 
object of \C, while each coupon of \GammaM\ is
labeled by an element of the morphism space of \C\ that corresponds
to the objects of the ribbons which enter and leave the coupon.
Composition in $\Cob_\C$ is defined by gluing, the identity
morphism $\id_\E$ is the cylinder over \E, and the tensor
product is given by disjoint union of objects and cobordisms.

A topological field theory furnishes, for any extended surface,
a representation of the mapping class group. Our approach relies on the
fundamental conjecture (which is largely established for a broad class of models)
that, for \C\ the representation
category of a rational vertex algebra \V, the mapping class group representation 
given by \tftC\ is equivalent to the one provided by the Knizhnik-Zamolodchikov 
connection on the conformal blocks for the vertex algebra \V.


\subsection{The TFT construction of full RCFT}

Let us now turn to the discussion of full local conformal field theories,
which are the structures to be compared to holonomies.
A {\em full\/} CFT is, by definition, a consistent system
of {\em local correlation functions\/} that satisfy all sewing constraints
(see e.g.\ \cite[def.\,3.14]{fjfrs2}). According to the principle of holomorphic 
factorization, every full RCFT can be understood with the help of a 
corresponding chiral CFT. The relevant chiral CFT is, however,
not defined on world sheets \X\ (which may be unoriented or have a non-empty
boundary), but rather on their {\em complex doubles\/} \Xh, which
can be given the structure of extended surfaces; 
this affords a geometric separation of left- and right-movers.
The double \Xh\ of \X\ is, by definition, 
the orientation bundle over \X\ mo\-du\-lo identification of
the two points in the fibre over each boundary point of \X. The world sheet
\X\ can be obtained from \Xh\ as the quotient by an orientation-reversing 
involution $\tau$. To give some examples, when \X\ is closed and orientable,
then \Xh\ is just the disconnected sum $\Xh\eq\X\,{\sqcup}\,{-}\X$ of 
two copies of \X\ with opposite orientation, and the involution $\tau$ 
just exchanges these two copies; the double of both the 
disk and the real projective plane is the two-sphere (with $\tau$ 
being given, in standard complex coordinates, by $z{\mapsto}\overline z^{-1}$
and by $z{\mapsto}\!{-}\overline z^{-1}$, respectively); and the double of
both the annulus and the M\"obius strip is a two-torus.
Further, when \X\ comes with field insertions, that is, embedded arcs
labeled by objects of either \C\ (for arcs on $\partial\X$) or pairs
of objects of \C\ (for arcs in the interior of \X), then corresponding arcs 
labeled by objects of \C\ are present on \Xh. 

\medskip

Given this connection between the surfaces relevant to chiral and full CFT, 
the relationship between the chiral and the full CFT can be stated as follows: 
A correlation function \CX\ of the full CFT on \X\ is a specific element in 
the appropriate space of conformal blocks of the chiral CFT on the double \Xh. 
A construction of such elements has been accomplished in 
\cite{fuRs4,fuRs8,fuRs10,fjfrs}. The first observation is that they
can be computed with the help of the corresponding TFT, namely as
  \be
  \CX = \tftC(\MX)\,1 ~\in \tftC(\Xh) \,.
  \label{fnsw:CX}\ee
Here $\MX\;{\equiv}\;\emptyset{\stackrel\MX\longrightarrow}\Xh$, 
the {\em connecting manifold\/} for the world sheet \X, is an extended 
cobordism that is constructed from the data of \X. Besides the category \C,
the specification of the vector \CX\ needs a second ingredient: a (Morita 
class of a) symmetric special Frobenius algebra $A$ in \C.
\\[3pt]
Let us give some details\,%
  \footnote{~For another brief summary, with different emphasis, see Section 7
  of \cite{fuRs11}. An in-depth exposition, including for instance the relevance
  of various orientations, can e.g.\ be found in Appendix B of \cite{fjfrs}. }
of the construction of \CX.
\def\leftmargini{1.12em}~\\[-1.6em]
\begin{itemize}
\itx
As a three-manifold, \MX\ is the interval bundle over \X\ modulo a 
$\mathbb{Z}_2$-identification of the intervals over $\partial\X$. 
Explicitly,
  \be
  \MX = \big( \Xh\Times[-1,1] \big) /{\sim}  \qquad{\rm with}\qquad
  ([x,\ort],t)\sim([x,-\ort],-t) \,.
  \ee
It follows in particular that $\partial\MX\eq\Xh$ and that \X\ is naturally 
embedded in \MX\ as $\imath\colon \X{\stackrel\simeq\to}\X\Times\{t{=}0\}\,
{\hookrightarrow}\MX$. Indeed, $\imath(\X)$ is a deformation retract of
\MX, so that the topology of \MX\ is completely determined by the one of \X.
\itx
A crucial ingredient of the construction of the ribbon graph \GammaMX\ in \MX\ is a 
(dual) oriented triangulation $\Gamma$ of the submanifold $\imath(\X)$ of \MX. 
This triangulation is labeled by objects and morphisms of \C. It is here that
the Frobenius algebra $A$ enters: Each edge of $\Gamma\,{\setminus}\,
\imath(\partial\X)$ is covered with a ribbon labeled by the object $A$ of \C, 
while each (three-valent) vertex is covered with a coupon labeled by the
multiplication morphism $m\iN\Hom_\C(A\oti A,A)$. In addition, whenever 
these assignments in themselves would be in conflict with the orientations 
of the edges, a coupon with morphism in either $\Hom_\C(A\oti A,\one)$ or
$\Hom_\C(\one,A\oti A)$ is inserted. Such morphisms are part of the data
for a Frobenius structure on $A$.
Assuming, for now, that the world sheet \X\ is oriented, 
independence of \CX\ from the choice of triangulation $\Gamma$ amounts
precisely to the statement that the object $A$ carries the structure
of a symmetric special Frobenius algebra.
\itx
If \X\ has {\em non-empty boundary\/}, the prescription for $\Gamma$ is
amended as follows. Each edge $e$ of $\Gamma\,{\cap}\,\imath(\partial\X)$ 
is covered with a ribbon labeled by a (left, say) $A$-module $N \eq N(e)$,
while each vertex lying on $\imath(\partial\X)$ is covered with a coupon 
that has incoming $N$- and $A$-ribbons as well as an outgoing $N$-ribbon and 
that is labeled by the representation morphism $\rho_N \iN\Hom_\C(A\oti N,N)$. 
The physical interpretation of the $A$-module $N$ is as the 
{\em boundary condition\/} that is associated to a component of 
$\partial\X$. That the object $N$ of \C\ labeling a boundary condition
carries the structure of an $A$-module and that the morphism $\rho_N$ is the 
corresponding representation morphism is precisely what is required (in 
addition to $A$ being a symmetric special Frobenius algebra) in order to get
independence of \CX\ from the choice of triangulation $\Gamma$.
\itx
If \X\ is {\em unoriented\/}, then as an additional feature one must 
ensure independence of \CX\ from the choice of local orientations of \X. As 
shown in \cite{fuRs8}, this requires an additional structure on the algebra $A$,
namely the existence of a morphism $\sigma \iN \Hom_\C(A,A)$ that
is an algebra isomorphism from the opposite algebra $A^{\rm opp}$ to
$A$ and squares to the twist of $A$, i.e.\ satisfies
  \be
  \sigma \circ \eta = \eta \,, \qquad
  \sigma \circ m = m \circ c_{A,A}^{} \circ (\sigma\oti\sigma) \,, \qquad
  \sigma \circ \sigma = \theta_{\!A} \,,
  \label{fnsw:sig}\ee
where $\eta\iN\Hom_\C(\one,A)$, $\theta_{\!A} \iN \Hom_\C(A,A)$ and
$c_{A,A} \iN \Hom_\C(A\oti A,A\oti A)$
denote the unit morphism, the twist, and the self-braiding of $A$, respectively.
This way $A$ becomes a braided version of an algebra with involution.
A symmetric special Frobenius algebra endowed with a morphism $\sigma$
satisfying \erf{fnsw:sig} is called a {\em Jandl algebra\/}.

\itx
In the presence of {\em topological defect lines\/} on \X\ a further amendment
of the prescription is in order. The defect lines partition \X\ into disjoint 
regions, and to the regions to the left and to the right of a defect line one may
associate different (symmetric special Frobenius) algebras $A_l$ and $A_r$,
such that the part of the triangulation $\Gamma$ in one region is labeled
by the algebra $A_l$, while the part of $\Gamma$ in the other region is labeled
by $A_r$. The defect lines are to be regarded as forming a subset
$\Gamma^{D}$ of $\Gamma$ themselves; each edge of $\Gamma^{D}$ 
is covered with a ribbon labeled by some object $B$ of \C, while each vertex of
$\Gamma$ lying on $\Gamma^{D}$ is covered with a coupon labeled by a morphism
$\rho\iN\Hom_\C(A_l\oti B,B)$, respectively $\rhor\iN\Hom_\C(B\oti A_r,B)$.
Consistency requires that these morphisms endow the object $B$ of \C\ that 
labels a defect line with the structure of an $A_l$-$A_r$-bimodule.
(Below we will concentrate on the case $A_l\eq A_r\,{=:}\,A$, so that we deal
with $A$-bimodules.)
\itx 
There are also rules for the morphisms of \C\ that label bulk, 
boundary and defect fields, respectively. 
\end{itemize}

The prescription summarized above allows one to construct 
the correlator \erf{fnsw:CX} for any arbitrary world sheet \X.
The so obtained correlators can be proven \cite{fjfrs} to satisfy all
consistency conditions that the correlators of a CFT must obey. Thus, 
specifying the algebra $A$ is sufficient to obtain a consistent system of 
correlators. The assignment of a (suitably normalized) correlator \CX\ to \X\ 
actually depends only on the Morita class of the symmetric special Frobenius 
algebra $A$. Conversely, any consistent set of correlators can be obtained
this way \cite{fjfrs2}.

Topological defects admit a number of interesting operations. In particular,
they can be fused -- on the algebraic side this corresponds to the tensor
product $B\,{\otimes_A}\, B'$ of bimodules. The bimodule morphisms 
$\Hom_{A|A}(B\,{\otimes_A}\,B',B'')$ appear as labels of vertices of defect 
lines. Defect lines can also be fused to boundaries; depending on the relative 
situation of the defect line and the boundary, this is given on the algebraic 
side by the tensor product $B\,{\otimes_A}\,N$ of a bimodule with a left module, 
or by the tensor product $N\,{\otimes_A}\,B$ with a right module, respectively.

\medskip

In the following table we collect some pertinent aspects of the construction and 
exhibit the geometric structures on the sigma model target space \Mt\ that 
correspond to them. 
\\[-2.6em]

\begin{center}
\begin{tabular}{l|l|l}
~geometric situation & algebraic structure in the category \C & \hspace*{1.1em}geometric 
structure on \Mt 
\\[-.9em]&\\\hline&\\[-.8em]
\hspace*{-.4em}\X\ closed oriented & symm.\ special Frobenius algebra $A$
                   & bundle gerbe $\mathcal{G}$ with connection \hspace*{-1em}
\\[2pt]
\hspace*{-.4em}\X\ unoriented & Jandl structure $\,\sigma\colon A^{\rm opp} \To A$
                   & Jandl gerbe 
\\[2pt]
\hspace*{-.4em}boundary condition & $A$-module & $\mathcal{G}$-D-brane 
\\[2pt]
\hspace*{-.4em}topological defect line & $A$-bimodule & $\mathcal{G}$-bi-brane
\end{tabular}
\end{center}

\noindent
Jandl gerbes, D-branes and bi-branes will be
presented in Sections \ref{fnsw:sec5}, \ref{fnsw:sec6} and \ref{fnsw:sec7}, respectively.


\section{Jandl gerbes: Holonomy for unoriented surfaces} \label{fnsw:sec5}

We have defined trivial bundle gerbes with connection as 2-forms because 2-forms 
can be integrated over oriented surfaces. Closing the 2-category of trivial 
bundle gerbes under descent has lead us to bundle gerbes. {\em Jandl gerbes\/} are 
bundle gerbes with additional structure, whose holonomy is defined 
for closed surfaces without orientation, even for unorientable surfaces \cite{scsW}. 
In particular, Jandl gerbes provide Wess-Zumino terms for unoriented surfaces. Comparing 
the geometric data with the representation theoretic ones from Section \ref{fnsw:sec4}, 
bundle gerbes with connection correspond to Frobenius algebras, while Jandl gerbes 
correspond to Jandl algebras.

The appropriate quantity that has to replace 2-forms in order to make integrals 
over an unoriented surface well-defined is a \emph{2-density}. Every surface 
$\Sigma$ has an oriented double covering $\mathrm{pr}\colon \hat\Sigma \To \Sigma$
that comes with an orientation-reversing involution $\sigma\colon \hat\Sigma \To 
\hat\Sigma$ which exchanges the two sheets and preserves the fibres. A 2-density 
on $\Sigma$ is a 2-form $\omega\iN\Omega^2(\hat\Sigma)$ such that
  \begin{equation}
  \label{fnsw:1}
  \sigma^{*}\omega =- \omega\,\text{.}
  \end{equation}
Every ordinary 2-form $\rho$ on $\Sigma$ defines a particular 2-density by
$\omega_{\rho} \,{:=}\, \mathrm{pr}^{*}\rho$.

\medskip

A 2-density on $\Sigma$ can indeed be integrated without requiring $\Sigma$ to be 
oriented. One chooses a dual triangulation $\Gamma$ of $\Sigma$ and, for each face $f$ 
of $\Gamma$, one of its two preimages under $\mathrm{pr}\colon\hat\Sigma \To \Sigma$, 
denoted $f_{\mathrm{or}}$. Then one sets
  \begin{equation}
  \label{fnsw:2}
  \int_{\Sigma} \omega := \sum_{f} \int_{\!f_{\mathrm{or}}} \omega \,\text{.}
  \end{equation}
Owing to the equality \erf{fnsw:1} the so defined integral does not depend on the 
choice of the preimages $f_{\mathrm{or}}$ nor on the choice of triangulation 
$\Gamma$. If $\Sigma$ can be endowed with an orientation,
the preimages $f_{\mathrm{or}}$ can be chosen 
in such a way that $\mathrm{pr}|_{f_{\mathrm{or}}}\colon f_{\mathrm{or}} \To f$ 
is orientation-preserving. Then the integral of a 2-density $\omega_{\rho} $ 
coincides with the ordinary integral of the 2-form $\rho$. 

Next we want to set up a 2-category whose objects are related to 2-densities. 
To this end we use the 2-category of 
trivial bundle gerbes introduced in Section \ref{fnsw:sec:defgrb}. 
Thus, one datum specifying an object is a 2-form $\omega \iN \Omega^2(\hat\Sigma)$.
In the context of 2-categories, demanding strict equality as in \erf{fnsw:1}
is unnatural. Instead, we replace equality by a 1-morphism
  \begin{equation}
  \label{eta}
  \eta:\quad \sigma^{*}\omega \,\to {-}\omega\,\text{,}
  \end{equation}
i.e.\ a 1-form $\eta\iN\Omega^1(\hat\Sigma)$ such that $\sigma^{*}\omega \eq {-} 
\omega + \mathrm{d}\eta$. As we shall see in a moment, we must impose equivariance 
of the 1-morphism up to some 2-morphism, i.e.\ we need in addition a 2-isomorphism
  \begin{equation}
  \label{phi}
  \phi:\quad \sigma^{*}\eta {\Rightarrow} \eta\,\text{,}
  \end{equation}
in other words a smooth function $\phi\colon M \To \Ue$ such that $\eta \eq 
\sigma^{*}\eta {-} \ii\,\mathrm{dlog}\,\phi$. This 2-iso\-mor\-phism, in turn, must 
satisfy the equivariance relation
  \begin{equation}
  \label{fnsw:3}
  \sigma^{*}\phi=\phi^{-1}\text{.}
  \end{equation}

Thus the objects of the 2-category are triples $(\omega,\eta,\phi)$. Let us verify 
that they still lead to a well-defined notion of holonomy. We choose again a dual 
triangulation $\Gamma$ of $\Sigma$ as well as a preimage $f_{\mathrm{or}}$ 
for each of its faces. The expression \erf{fnsw:2} is now no longer 
independent of these choices, because every change creates a boundary term in 
the integrals of the 1-form $\eta$. To resolve this problem, we
involve \emph{orientation-reserving edges}: these are edges in $\Gamma$
whose adjacent faces have been lifted to opposite sheets. Since 
$\Gamma$ is a dual triangulation, its orienta\-ti\-on-re\-ver\-sing edges form a 
disjoint union of piecewise smooth circles $c\,{\subset}\,\Sigma$. For each of these 
circles, we choose again a preimage $c_{\mathrm{or}}$. It may not be possible to 
choose $c_{\mathrm{or}}$ to be closed, in which case
there exists a point $p^c\iN\Sigma$ which has two preimages in $c_{\mathrm{or}}$.
We choose again one of these preimages, denoted $p^c_{\mathrm{or}}$. Then
  \begin{equation}
  \label{fnsw:4}
  \hol_{\omega,\eta,\phi} := \exp\! \left( 2\pi\im \Big (
  \sum_{f} \int_{\!f_{\mathrm{or}}} \omega + \sum_{c} \int_{c_{\mathrm{or}}} \eta \Big) \right)
  \prod_{c} \phi(p^{c}_{\mathrm{or}})
  \end{equation}
is independent of the choice of the lifts $f_{\mathrm{or}}$, $c_{\mathrm{or}}$ 
and $p_{\mathrm{or}}$, and is independent of the choice of the triangulation. 

More generally, let $\MI$ be the category of smooth manifolds with involution, whose 
morphisms are equivariant smooth maps. (The involution is not required to act freely.)
In a first step, we want to define a presheaf
  \begin{equation}
  \Jantriv :\quad \MoppI \to \Cat
  \end{equation}
of trivial Jandl gerbes. For $(M,k)$ a smooth manifold with involution 
$k\colon M \To M$, a trivial Jandl gerbe involves as a first datum
a trivial bundle gerbe $\mathcal{I}_{\omega}$, but  
as explained in Section \ref{fnsw:def:grb} we replace the 1-morphism $\eta$ from 
\erf{eta} by a line bundle $L$ over $M$ with connection of curvature 
  \begin{equation}
  \label{fnsw:5}
  \mathrm{curv}(L) = -\,\omega - k^{*}\omega\,\text{,}
  \end{equation}
and we replace the 2-isomorphism $\phi$ from \erf{phi} by an isomorphism 
$\phi\colon k^{*}L \To L$ of 
line bundles with connection, still subject to the condition \erf{fnsw:3}. 
Notice that the pair $(L,\phi)$ is nothing but a
$k$-equivariant line bundle with connection over $M$. After this  step, we 
still have the holonomy \erf{fnsw:4}, which now looks like
  \begin{equation}
  \hol_{\mathcal{I}_{\omega},L,\phi}
  = \exp \!\Big ( 2\pi \im \sum_{f}\int_{\!f_{\mathrm{or}}}\omega \Big )\,
  \prod_c \hol_{\bar L}(c)\,\text{,}
  \end{equation}
where we have used the fact that, since the action of $\langle k\rangle$ on 
$c_{\mathrm{or}}$ is free, the $k$-equivariant line bundle $(L,\phi)$ descends 
to a line bundle $\bar L$ with connection over the quotient 
$c \eq c_{\mathrm{or}}/\langle k\rangle$. 
This formula is now manifestly independent of the 
choices of $c_{\mathrm{or}}$ and $p^c_{\mathrm{or}}$. Its independence under 
different choices of faces $f_{\mathrm{or}}$ is due to \erf{fnsw:5}. 

Now we close the presheaf $\Jantriv(M)$ under descent to allow for
non-trivial bundle gerbes. To do so, we need to introduce duals of bundle 
gerbes, 1-morphisms and 2-isomorphisms see \cite{waldorf1}; for the sake of 
brevity we omit these definitions here.

\begin{definition}
Let $M$ be a smooth manifold with involution $k\colon M \To M$. A \emph{Jandl gerbe} 
is a bundle gerbe $\mathcal{G}$ over $M$ together with a 1-isomorphism 
$\mathcal{A}\colon k^{*}\mathcal{G} \To \mathcal{G}^{*}$ to the dual gerbe 
and a 2-isomorphism $\varphi\colon k^{*}\mathcal{A}{\Rightarrow}\mathcal{A}^{*}$ 
that satisfies $k^{*}\varphi \eq \varphi^{*-1}$.
\end{definition}

Jandl gerbes form a sheaf
  \begin{equation}
  \Jan:\quad \MoppI \to \Cat\,\text{.}
  \end{equation}
The gluing axiom for this sheaf has been proved in \cite{gasuwa2}. We remark that the 
1-isomorphism $\mathcal{A}$ may be regarded as the counterpart of a Jandl structure 
$\sigma$ on the Frobenius algebra $A$ that corresponds to the bundle gerbe $\mathcal{G}$, 
if one accepts 
that the dual gerbe plays the role of the opposed algebra.       

\medskip

Suppose we are given a Jandl gerbe $\mathcal{J}$ over a smooth manifold $M$ with 
involution $k$. 
If $\Sigma$ is a closed surface, possibly unoriented and possibly unorientable, and
  \begin{equation}
  \Phi:\quad (\hat\Sigma,\sigma) \to (M,k)
  \end{equation}
is a morphism in $\MoppI$, we can pull back the Jandl gerbe $\mathcal{J}$ from $M$ to 
$\hat\Sigma$. As in the case of ordinary surface holonomy, it then becomes trivial 
as a gerbe for dimensional reasons, and we can choose an isomorphism
  \begin{equation}
  \mathcal{T}:\quad \Phi^{*}\mathcal{J} \iso (\mathcal{I}_{\omega},L,\phi)\,\text{.}
  \end{equation}
Then we define
  \begin{equation}
  \hol_{\mathcal{J}}(\Phi) := \hol_{\mathcal{I}_{\omega},L,\phi}\,\text{.}
  \end{equation}
This is independent of the choice of $\mathcal{T}$, because any other choice 
$\mathcal{T}'$ gives rise to an isomorphism $\mathcal{T}'\circ\mathcal{T}^{-1}$ in 
$\Jantriv(\hat\Sigma,\sigma)$ under which the holonomy stays unchanged.

\medskip

We have now seen that every Jandl gerbe $\mathcal{J}$ over a smooth manifold $M$ 
with involution $k$ has holonomies for unoriented closed surfaces and equivariant 
smooth maps $\Phi\colon \hat\Sigma \To M$. We thus infer that  sigma models on $M$ 
whose fields are such maps, are defined by Jandl gerbes $\mathcal{J}$ over $M$ 
rather than by ordinary bundle gerbes $\mathcal{G}$. This makes it an interesting 
problem to classify Jandl gerbes.  

Concerning the existence of a Jandl gerbe $\mathcal{J}$ with underlying bundle gerbe 
$\mathcal{G}$, the 1-iso\-morphism $\mathcal{A}\colon k^{*}\mathcal{G} \To 
\mathcal{G}^{*}$ requires the curvature $H$ of $\mathcal{G}$ to satisfy
  \begin{equation}
  \label{fnsw:6}
  k^{*}H = -H\,\text{.}
  \end{equation}
Apart from this necessary condition, there is a sequence of obstruction classes 
\cite{gasuwa2}. Reduced to the case that $M$ is 2-connected, there is one obstruction 
class $o(\mathcal{G}) \iN H^3(\Z_2,\Ue)$, the group cohomology of $\Z_2$ with 
coefficients in $\Ue$, on which $\Z_2$ acts by inversion. If $o(\mathcal{G})$ vanishes,
then inequivalent Jandl gerbes with the same underlying bundle gerbe $\mathcal{G}$ 
are parameterized by $H^2(\Z_2,\Ue)$. 

These results can be made very explicit in the case of WZW models, 
for which the object in $\MI$ is a connected compact simple Lie group $G$ 
equipped with an involution $k\colon G\To G$ acting as
  \begin{equation}
  k\colon\quad g \,\mapsto (zg)^{-1}
  \end{equation}
for a fixed `twist element' $z\iN Z(G)$. It is easy to see that the  3-form 
$H_k \iN\Omega^3(G)$, which is the curvature of the level-$k$ bundle gerbes 
$\mathcal{G}$ over $G$, satisfies the necessary condition \erf{fnsw:6}. All obstruction 
classes $o(\mathcal{G})$ and all parameterizing groups have been computed in 
dependence of the twist element $z$ and the level $k$ \cite{gasuwa1}. The 
numbers of inequivalent Jandl gerbes  range from two (for simply connected $G$,
per level and involution) to sixteen (for $\mathrm{PSO}(4n)$, for every even level). 

Most prominently, there are two involutions on $\mathrm{SU}(2)$, namely 
$g{\mapsto}g^{-1}$ and $g{\mapsto}{-}g^{-1}$, and for each of them two inequivalent 
Jandl gerbes per level. On $\mathrm{SO}(3)$ there is only a single involution, but 
the results of \cite{scsW,gasuwa1} exhibit four inequivalent Jandl gerbes per 
even level. This explains very nicely why $\mathrm{SU}(2)$ and $\mathrm{SO}(3)$ 
have the same number of orientifolds, despite a different 
number of involutions. These results reproduce those of the algebraic approach
(see e.g.\ \cite{fuRs8}); for the precise comparison, Jandl structures related 
by the action of the trivial line bundle with either of its two equivariant 
structures have to be identified.


\section{D-branes: Holonomy for surfaces with boundary} \label{fnsw:sec6}

We now introduce the geometric structure needed to define surface holonomies and 
Wess-Zumino terms for surfaces with boundary. When one wants to define holonomy 
along a curve that is not closed, one way to make the parallel transport group-valued 
is to choose trivializations at the end points. 
To incorporate these trivializations into the background,
one can choose a submanifold $\dQ \,{\subset}\, M$ together with a 
trivialization $E|_\dQ \To \trivlin_A$. Admissible paths $\gamma\colon [0,1] \To M$ are 
then required to start and end on this submanifold, $\gamma(0), \gamma(1) \iN \dQ$. 
The same strategy has proven to be successful for \emph{surfaces} with boundary.

\begin{definition}
Let $\mathcal{G}$ be a bundle gerbe with connection over $M$.
A \emph{$\mathcal{G}$-D-brane} 
is a submanifold $\dQ \,{\subset}\, M$ together with a 1-morphism 
  \be
  \mathcal{D}\colon\quad \mathcal{G}|_\dQ \to \mathcal{I}_{\omega}
  \ee
to a trivial bundle gerbe $\mathcal{I}_{\omega}$ given by a two-form $\omega$ on $\dQ$.
\end{definition}

The morphism $\mathcal{D}$ is called a $\mathcal{G}$-module, or 
twisted vector bundle. Notice that if $H$ is the curvature of $\mathcal{G}$, 
the 1-morphism $\mathcal{D}$ enforces the identity
  \begin{equation}
  H|_\dQ = \mathrm{d}\omega\,\text{.}
  \label{fnsw:H=domega}\end{equation}
This equality restricts the possible choices
of the world volume $\dQ$ of the $\mathcal{G}$-D-brane.

Suppose that $\Sigma$ is an oriented surface, possibly with boundary, and 
$\Phi\colon \Sigma \To M$ is a smooth map. We require that 
$\Phi(\partial\Sigma) \,{\subset}\, \dQ$. As described in Section 
\ref{fnsw:subsec:surfhol}, we choose a trivialization $\mathcal{T}\colon 
\Phi^{*}\mathcal{G} \To \mathcal{I}_{\rho}$. Its restriction to 
$\partial\Sigma$
and the $\mathcal{G}$-module $\mathcal{D}$ define a 1-morphism
  \begin{equation}
  \alxydim{} {\mathcal{I}_{\rho}\big|_{\partial\Sigma}^{}
  \ar[rr]^{\mathcal{T}^{-1\!}|_{\partial\Sigma}^{}\hspace*{3em}} && \,
  \Phi^*\mathcal{G}\big|_{\partial\Sigma}^{} = \Phi^*(\mathcal{G}\big|_{\dQ}^{})
  \ar[rr]^{\hspace*{2.1em}\Phi^*(\mathcal{D})} && \Phi^*(\mathcal{I}_{\omega})\,\text{.}}
  \end{equation}
According to the definition \erf{defhomvect}, this 1-morphism is nothing but a 
hermitian vector bundle $E$ with connection over $\partial\Sigma$
and its curvature is $\mathrm{curv}(E) \eq \omega - \rho$. Then we consider
  \begin{equation}
  \mathrm{Hol}_{\mathcal{G},\mathcal{D}}(\Phi)
  := \exp \!\Big( 2\pi\im \int_{\Sigma} \rho \Big) \,
  \mathrm{tr}(\mathrm{Hol}_E(\partial\Sigma))\,\text{,}
  \end{equation}
where the trace makes the holonomy of $E$ independent of the choice of a parameterization 
of $\partial\Sigma$. This expression is independent of the choice of the trivialization 
$\mathcal{T}$: if $\mathcal{T}'\!\colon \mathcal{G} \To \mathcal{I}_{\rho'}$ is another 
one and $E'$ is the corresponding vector bundle, we have the transition isomorphism 
$L$ from \erf{fnsw:transiso} with curvature $\rho'-\rho$, and an isomorphism 
$E' \oti L \cong E$. It follows that
  \begin{equation}
  \exp \!\Big( 2\pi\im \int_{\Sigma} \rho \Big)
  \, \mathrm{tr}(\mathrm{Hol}_{E}(\partial \Sigma))
  = \exp \!\left (2\pi\im \Big ( \int_{\Sigma} \rho' - \mathrm{curv}(L) \Big) \right )
  \, \mathrm{tr}(\mathrm{Hol}_{E' \otimes L}(\partial \Sigma)) \,,
  \end{equation}
and on the right hand side the unprimed quantities cancel by Stokes' theorem. 

\medskip

Important results on D-branes concern in particular two large classes of 
models, namely free field theories and again WZW theories.
The simplest example of a free field theory is the one of a compactified free boson,
in which $M$ is a circle $S_R^1 \cong \R\bmod2\pi R\,\Z$ of radius $R$. As is well 
known, there are then in particular
two distinct types of D-branes: D0-branes $\mathcal D^{(0)}_x$, 
whose support is localized at a position $x\iN S_R^1$, and D1-branes $\mathcal D^{(1)}_\alpha$, 
whose world volume is all of $S_R^1$ and which are characterized by a Wilson line
$\alpha\iN \R\bmod\!\frac1{2\pi R}\,\Z$, corresponding to a flat connection on $S^1_R$.

For WZW theories, which are governed by a bundle gerbe $\mathcal{G}$ over a connected 
compact simple Lie group $G$, preserving the non-abelian
current symmetries puts additional constraints on the admissible D-branes: 
their support $\dQ$ must be a conjugacy class $\mathcal{C}_{h}$ of a group element $h\iN G$.
This can e.g.\ be seen by studying the scattering of bulk fields in 
the presence of the D-brane. On such conjugacy classes one finds a canonical 
2-form $\omega_{h} \iN \Omega^2(\mathcal{C}_{h})$. 
Additionally, the 1-morphism $\mathcal{D}\colon \mathcal{G}|_{\mathcal{C}_h} \To 
\mathcal{I}_{\omega_h}$ of a symmetric D-brane must satisfy a certain equivariance 
condition \cite{gawedzki4}. Interestingly, only on those conjugacy classes $\mathcal{C}_{h}$ for which 
  \begin{equation}
  h = \exp(2\pi\im\, \Frac{\alpha+\rho}{k+g^{\vee}}) \,,
  \label{fnsw:h-alpha}
  \end{equation}
with $\alpha$ an \emph{integrable} highest weight, admit such 1-morphisms. Here $\rho$ denotes 
the Weyl vector and $g^\vee$ the dual Coxeter number of the Lie algebra $\g$ of 
$G$. Thus in particular the possible world volumes of symmetric D-branes 
form only a discrete subset of conjugacy classes.

\medskip

We finally remark that the concepts of D-branes and Jandl gerbes can be merged 
\cite{gasuwa2}. The resulting structures provide holonomies for unoriented surfaces 
with boundary, and can be used to define D-branes in WZW orientifold theories.


\section{Bi-branes: Holonomy for surfaces with defect lines} \label{fnsw:sec7}

\subsection{Gerbe bimodules and bi-branes}

In the representation theoretic approach to rational conformal field theory, 
boundary conditions and defect lines are described as modules and bimodules, 
respectively. The fact that the appropriate target space structure for 
describing boundary conditions, D-branes, is related to gerbe modules, raises 
the question of what the appropriate target space structure for defect lines 
should be. The following definition turns out to be appropriate.

\begin{definition}
Let $\mathcal{G}_1$ and $\mathcal{G}_2$ be bundle gerbes with connection 
over $M_1$ and $M_2$, respectively. A \emph{$\mathcal{G}_1$-$\mathcal{G}_2$-bi-brane}
is a submanifold $\bQ \,{\subset}\, M_1 \Times M_2$ together with a
$(p_1^{*}\mathcal{G}_1)|_\bQ$-$(p_2^{*}\mathcal{G}_2)|_\bQ$-bimodule, i.e.\
with a 1-morphism
  \be
  \mathcal{B} :\quad
  (p_1^{*}\mathcal{G}_1^{})|_\bQ^{} \to (p_2^{*}\mathcal{G}_2^{})|_\bQ^{}
  \otimes \mathcal{I}_{\bbf} 
  \label{fnsw:D}\ee
with $\mathcal{I}_{\bbf}$ a trivial bundle gerbe given by a two-form $\bbf$ on $\bQ$.
\end{definition}

Similarly as in \erf{fnsw:H=domega} it follows that the two-form $\bbf$ on $\bQ$ obeys
  \be
  p_1^{*}H|_\bQ^{} = p_2^{*}H|_\bQ^{} + \mathrm{d}\bbf \,\text{.}
  \label{fnsw:p1=p2+domega}
  \ee
We call $\bQ$ the world volume and $\bbf$ the
curvature of the bimodule. With the appropriate notion of duality for bundle gerbes
(see Section 1.4 of \cite{waldorf1}), a $\mathcal{G}_1$-$\mathcal{G}_2$-bimodule
is the same as a $(\mathcal{G}_1{\otimes}\mathcal{G}_2^{*})$-module.
For a formulation in terms of local data, see (B.8) of \cite{fusW}.

\medskip

As an illustration, consider again the free boson and WZW theories,
restricting attention to the case $M_1\eq M_2$. For the free boson
compactified on a circle $S_R^1$ of radius $R$, one finds that
the world volume of a bi-brane is a submanifold 
$\bQ_x \,{\subset}\,S_R^1\Times S_R^1$ of the form
  \be
  \bQ_{x,\alpha} := \{ (y,y{-}x) \,|\, y \iN S_R^1 \}
  \label{fnsw:bQ-x-alpha}\ee
with $x\iN S_R^1$. The submanifold $\bQ_{x,\alpha}$ has the topology of a circle and 
comes with a flat connection, i.e.\ with a Wilson line $\alpha$. Thus the bi-branes of a 
compactified free boson are naturally parameterized by a pair $(x,\alpha)$ taking
values in two dual circles that describe a point on $S^1_R$ and a Wilson line.

\medskip

In the WZW case, for which the target space is a compact connected simple Lie group $G$,
a scattering calculation \cite{fusW} similar to the one performed for D-branes indicates 
that the world volume of a (maximally symmetric) bi-brane is a {\em biconjugacy class\/}
  \be
  \bQ_{h,h'} := \left\lbrace (g,g')\iN G \Times G \,|\,
  \exists\, x_1^{},x_2^{}\iN G\text{: }g\eq x_1^{}h^{}x_2^{-1}\text{,}\,
  g'^{}\eq x_1^{} h'^{} x_2^{-1} \right\rbrace
  \,\subset G\Times G
  \ee
of a pair $(h,h')$ of group elements satisfying $h\,(h')^{-1} \iN \mathcal{C}_{h_\alpha}$
with $h_\alpha$ as given in \erf{fnsw:h-alpha}. The bi\-con\-jugacy classes carry two 
commuting $G$-actions, corresponding to the presence of two independent conserved 
currents in the field theory.  Further, a biconjugacy class can be described as the preimage
  \be
  \bQ_{h,h'}^{} = \tilde\mu^{-1}(\mathcal{C}_{hh'{}^{-1}}) =
  \big\lbrace (g,g') \iN G \Times G \,|\, gg'{}^{-1}
  \iN \mathcal{C}_{hh'{}^{-1}} \big\rbrace 
  \ee
of the conjugacy class $\mathcal{C}_{hh'{}^{-1}}$ under the map
  \begin{equation}
  \tilde \mu :\quad  G \Times G \ni\, (g_1,g_2) \mapsto g_1^{}g_2^{-1} \,\in G\,\text{.}
  \end{equation}
Finally, the relevant two-form on $\bQ_{h,h'}$ is
  \be
  \bbf_{h,h'}:= \tilde\mu^* \omega_{hh'{}^{-1}} -
  \Frac{k}{2}\,\langle p_1^*\theta\wedge p_2^*\theta\rangle \,.
  \label{fnsw:bbf-h-h'}
  \ee
Here $k$ is the level, $\theta$ is the left-invariant Maurer-Cartan form, $p_i$ 
are the projections to the factors of $G\Times G$, and $\omega_h$ is the canonical 
2-form (see Section \ref{fnsw:sec6}) on the conjugacy class $\mathcal{C}_h$. One 
checks that $\bbf_{h,h'}$ is bi-invariant and satisfies \erf{fnsw:p1=p2+domega}. 

Examples of symmetric bi-branes can be constructed from symmetric D-branes using a 
multiplicative structure on the bundle gerbe $\mathcal{G}$ \cite{waldorf2}.
Another important class of examples are Poincar\'e line bundles.
These describe T-dualities; an elementary relation between 
T-duality and Poincar\'e line bundles is
provided \cite{saSc} by the equation of motion \cite{ruSu} in the presence of defects.


\subsection{Holonomy and Wess-Zumino term for defects}

The notion of bi-brane allows one in particular to define holonomy also
for surfaces with defect lines.

The simplest world sheet geometry involving a defect line consists of
a closed oriented world sheet $\Sigma$ together with an embedded oriented circle 
$S\,{\subset}\, \Sigma$ that separates the world sheet into two 
components, $\Sigma \eq \Sigma_1 \,{\cup_{S}}\, \Sigma_2$. 
Assume that the defect $S$ separates regions that support conformally invariant
sigma models with target spaces $M_1$ and $M_2$, respectively, and consider maps
$\phi_i{:}\ \Sigma_i \To M_i$ for $i\iN\{1,2\}$ such that the image of 
  \be
  \begin{array}{lrcl}
  \phi_S:& S &\!\!\to\!\!& M_1 \Times M_2 \\[3pt]
  & s &\!\!\mapsto\!\!& (\phi_1(s) , \phi_2(s))
  \eear \ee
is contained in the submanifold $\bQ$ of $M_1 \Times M_2$.
The orientation of $\Sigma_i$ is the one inherited
from the orientation of $\Sigma$, and without loss of generality we take
$\partial \Sigma_1 \eq S$ and $\partial \Sigma_2 \eq {-}S$.

We wish to find the Wess-Zumino part of the sigma model action,
or rather the corresponding holonomy $\hol_{\mathcal{G}_1,\mathcal{G}_2,\mathcal{B}}$, 
that corresponds to having bundle gerbes
$\mathcal{G}_1$ and $\mathcal{G}_2$ over $M_1$ and $M_2$
and a $\mathcal{G}_1$-$\mathcal{G}_2$-bi-brane $\mathcal{B}$.
The pullback of the bimodule \erf{fnsw:D} along the map $\phi_S{:}\ S \To \bQ$ 
gives a $(\phi_1^{*}\mathcal{G}_1^{})|_S^{}$-$(\phi_2^{*}\mathcal{G}_2^{})|_S^{}$-bimodule
  \be
  \phi_S^{*}\mathcal{B}:\quad (\phi_1^{*}\mathcal{G}_1^{})|_S^{} \to
  (\phi_2^{*}\mathcal{G}_2^{})|_S^{} \otimes \mathcal{I}_{\phi_S^{*}\bbf}
  \,\text{.}
  \ee
The pullback bundle gerbes $\phi_i^{*}\mathcal{G}_i$ over $\Sigma_i$ are
trivializable for dimensional reasons, and a choice $\mathcal{T}_i{:}\
\phi_i^{*}\mathcal{G}_i^{} \,{\to}\, \mathcal{I}_{\rho}$ of trivializations
for two-forms $\rho_i$ on $\Sigma_i$ produces a vector bundle $E$ over $S$.
We then define
  \be
  \hol_{\mathcal{G}_1,\mathcal{G}_2,\mathcal{B}}^{}(\Sigma,S)
  := \exp \!\Big ( 2\pi\mathrm{i}\! \int_{\Sigma_1}\! \rho_1 \Big ) \, \exp
  \Big (2\pi \mathrm{i}\! \int_{\Sigma_2}\! \rho_2 \Big ) \,
  \mathrm{tr}(\hol_{E}^{}(S)) \,\in \complex 
  \ee
to be the holonomy in the presence of the bi-brane $\mathcal{B}$.
As shown in Appendix B.3 of \cite{fusW}, for similar reasons as in the case
of D-branes the number $\hol_{\mathcal{G}_1,\mathcal{G}_2,\mathcal{B}}(\Sigma,S)$
is independent of the choice of the trivializations
$\mathcal{T}_1$ and $\mathcal{T}_2$.


\subsection{Fusion of defects}

In the field theory context of section \ref{fnsw:sec4} there are
natural notions of the fusion of a defect (an $A$-bi\-module) with a boundary
condition (a left $A$-module), yielding another boundary condition, and of the 
fusion of two defects, yielding another defect. Both of these are provided by 
the tensor product over the relevant Frobenius algebra $A$. These representation
theoretic notions of fusion have a counterpart on the geometric side as well.

Consider first the fusion of a defect with a boundary condition.
We allow for the general situation of a defect described by an
$M_1$-$M_2$-bi-brane with different target spaces $M_1$ and $M_2$. Thus take an 
$M_1$-$M_2$-bi-brane with world volume $\bB\,{\subseteq}\, M_1\Times M_2$ and 
an $M_2$-D-brane with world volume $\dD \,{\subseteq}\, M_2$.
The action of correspondences on sheaves suggests the following ansatz for
the world volume of the fusion product:
  \be
  (\mathcal B\star \mathcal D)\dot{} := p_1^{}\!\left( \bB \,{\cap}\, p_2^{-1}(\dD) \right)
  \label{fnsw:D-bi-fusion}\ee
with $p_i$ the projection $M_1\Times M_2\to M_i$. 
The corresponding ansatz for the fusion of an $M_1$-$M_2$-bi-brane $\mathcal B$ of 
world volume $\bB$ with an $M_2$-$M_3$-bi-brane $\mathcal B'$ of world volume $\bB'$ 
uses projections $p_{ij}$ from $M_1\Times M_2\Times M_3$ to $M_i\Times M_j$:
  \be
  (\mathcal B \star \mathcal B')\dot{}
  := p_{13}^{}\big( p_{12}^{-1}(\bB)\,{\cap}\,p_{23}^{-1}(\bB') \big) \,.
  \label{fnsw:bi-bi-fusion}\ee

\medskip

In general one obtains this way only subsets, rather than 
submanifolds, of $M_1$ and $M_1\Times M_3$, respectively.
On a heuristic level one would, however, expect that owing to quantization 
of the branes a finite superposition of branes is selected, which
should then reproduce the results obtained in the field theory setting.

We illustrate this again with the two classes of models
already considered, i.e.\ free bosons and WZW theories, again
restricting attention to the case $M_1\eq M_2$.
First, for the theory of a compactified free boson, the D-brane 
is of one of the types $\mathcal D^{(0)}_x$ or $\mathcal D^{(1)}_\alpha$
(see Section \ref{fnsw:sec6}) and the bi-brane world volume is of the form
$\bQ_{x,\alpha}$ given in \erf{fnsw:bQ-x-alpha}. For D-branes of type 
$\mathcal D^{(0)}_x$ the prescription \erf{fnsw:D-bi-fusion} thus yields
  \be
  \mathcal B_{(x,\alpha)} \star \mathcal D^{(0)}_y = \mathcal D^{(0)}_{x+y} \,.
  \ee
For the fusion of a bi-brane $\mathcal B_{(x,\alpha)}$ and a D1-brane 
$\mathcal D^{(1)}_\beta$, one must
take the flat line bundle on the bi-brane into account. We first pull
back the line bundle on $\dQ^{(1)}_\beta$ along $p_2$ to a line
bundle on $S^1_R\Times S^1_R$, then restrict it to 
$\bQ_{(x,\alpha)}$, and finally tensor this restriction
with the line bundle on $\bQ_{(x,\alpha)}$ described by the Wilson
line $\alpha$. This results in a line bundle with Wilson line $\alpha{+}\beta$
on the bi-brane world volume, which in turn can be pushed down along 
$p_1$ to a line bundle on $S^1_R$, so that
  \be
  \mathcal B_{(x,\alpha)} \star \mathcal D^{(1)}_\beta = \mathcal D^{(1)}_{\alpha+\beta} \,.
  \ee
In short, the fusion with a defect $\mathcal B_{(x,\alpha)}$ acts on D0-branes
as a translation by $x$ in position space, and
on D1-branes as a translation by $\alpha$ in the space of Wilson lines.
Similarly, the prescription \erf{fnsw:bi-bi-fusion} leads to
  \be
  \mathcal B_{(x,\alpha)} \star \mathcal B_{(x',\alpha')}
  = \mathcal B_{(x+x',\alpha+\alpha')}
  \ee
for the fusion of two bi-branes $\mathcal B_{(x,\alpha)}$ and $\mathcal B_{(x',\alpha')}$,
i.e.\ both the position and the Wilson line variable of the bi-branes add up. 

\medskip

For WZW theories, besides the quantization of the positions of the branes
another new phenomenon is that multiplicities other than zero or one appear
in the field theory approach. In that context they
arise from the decomposition $B_\alpha \,{\otimes_{\!A}}\, B_\beta \eq
\bigoplus_\gamma\caln_{\alpha\beta}^{~~\gamma}B_\gamma$ of a tensor product
of simple $A$-bimodules into a finite direct sum of simple $A$-bimodules,
and analogously for the case of mixed fusion (in rational CFT, both the
category of $A$-modules and the category of $A$-bimodules are semisimple).
Moreover, for simply connected groups, the multiplicities appearing in both 
types of fusion are in fact the same as the chiral
fusion multiplicities which are given by the Verlinde formula.

By analogy with the field theory situation we expect fusion rules 
  \be
  \mathcal B_\alpha \,{\star}\, \mathcal B_\beta
  = \sum_\gamma \caln_{\alpha\beta}^{~~\gamma}\mathcal B_\gamma
  \label{fnsw:B*B}\ee
of bi-branes, and analogously for mixed fusion of bi-branes and D-branes.
In the particular case of WZW theories on simply connected Lie groups
one can in addition invoke the duality $\alpha{\mapsto}\alpha^\vee$ which 
in that 
case exists on the sets of branes as well as defects that preserve all current 
symmetries, so as to work instead with fusion coefficients of type
$\caln_{\alpha\beta\gamma}^{}\eq\caln_{\alpha\beta}^{}{}^{\!\gamma^\vee}$.
Then for the case of two D-branes $\mathcal D_\alpha$ and $\mathcal D_\gamma$
with world volumes given by conjugacy classes $\calc_{h_\alpha}$
and $\calc_{h_\gamma}$ of $G$, as well as a bi-brane
$\mathcal B_\beta$ whose world volume is the biconjugacy class
$\tilde\mu^{-1}(\mathcal{C}_{h_\beta})$, one is lead to consider the subset
  \be 
  \Pi_{\alpha\beta\gamma} 
  := p_1^{-1}(\calc_\alpha) \cap \tilde\mu^{-1}(\calc_\beta) \cap p_2^{-1}(\calc_\gamma)
  = \{ (g,g')\iN G\Times G \,|\, g\iN\calc_\alpha,\,
  g'\iN\calc_\gamma,\, gg'{}^{-1}{\in}\,\calc_\beta \}
  \ee
of $G\Times G$. Combining the adjoint action on $g$ and on $g'$ gives a natural $G$-action
on $\Pi_{\alpha\beta\gamma}$. And since both D-branes and 
the bi-brane are equipped with two-forms $\omega_\alpha$, $\omega_\gamma$ and 
$\bbf_\beta$, $\Pi_{\alpha\beta\gamma}$ comes with a natural two-form as well, namely with
  \be
  \omega_{\alpha\beta\gamma}
  := p_1^*\omega_\alpha|^{}_{\Pi_{\alpha\beta\gamma}}
  + p_2^*\omega_\gamma|^{}_{\Pi_{\alpha\beta\gamma}}
  + \bbf_\beta|^{}_{\Pi_{\alpha\beta\gamma}} \,.
  \label{fnsw:omegaabc}\ee

By comparison with the field theory approach, this result should
be linked to the fusion rules of the chiral WZW theory 
and thereby provide a physically motivated realization of the Verlinde algebra.
To see how such a relation can exist, notice that fusion rules are dimensions of
spaces of conformal blocks and as such can be obtained by geometric
quantization from suitable moduli spaces of flat connections which
arise in the quantization of Chern-Simons theories (see e.g.\ \cite{axdw}). The
moduli space relevant to us is the one for the three-punctured sphere $S^2_{(3)}$,
for which the monodromy of the flat connection around the punctures takes values 
in conjugacy classes $\calc_\alpha$, $\calc_\beta$ and $\calc_\gamma$, 
respectively. The relations in the fundamental group of $S^2_{(3)}$ imply
the condition $g_\alpha g_\beta^{} g_\gamma \eq 1$ on the
monodromies $g_\alpha\iN \calc_\alpha$, $g_\beta\iN\calc_\beta$ and
$g_\gamma\iN \calc_\gamma$. Since monodromies are defined only
up to simultaneous conjugation, the moduli space that matters in classical
Chern-Simons theory is isomorphic to the quotient $\Pi_{\alpha\beta\gamma}/G$.

It turns out that the range of bi-branes appearing in the fusion 
product is correctly bounded already before geometric quantization. Indeed,
the relevant product of conjugacy classes is
  \be
  \calc_h * \calc_{h'} := \{ gg' \,|\, g\iN\calc_h,\,g'\iN\calc_{h'} \} \,,
  \ee
and for the case of $G\eq$SU(2) it is easy to see that this yields
the correct upper and lower bounds for the SU(2) fusion rules \cite{jewe,fusW}.
A full understanding of fusion can, however, only be expected after applying geometric
quantization to the so obtained moduli space: this space must be endowed with
a two-form, which is interpreted as the curvature of a line bundle, and the
holomorphic sections of this bundle are what results from geometric
quantization. In view of this need for quantization it is a highly non-trivial
observation that the two-form \erf{fnsw:omegaabc} furnished by the two branes and
the bi-brane is exactly the same as the one that arises
from classical Chern-Simons theory.

\medskip

In terms of defect lines, the decomposition \erf{fnsw:B*B} of
the fusion product of bi-branes corresponds to the presence of a 
\emph{defect junction}, which constitutes a particular type of \emph{defect field}.
A sigma model description for world sheets with such embedded 
defect junctions has been proposed in \cite{ruSu}.

\medskip

We have demonstrated how structural analogies between the geometry of bundle 
gerbes and the representation theoretic approach to rational conformal field 
theory lead to interesting geometric structure, including a physically
motivated realization of the Verlinde algebra. The precise form of the latter
and its relation with the realization of the Verlinde algebra in the context 
of supersymmetric conformal field theory \cite{frht5} remain to be understood. 
But in any case the parallelism between classical actions and full quantum 
theory exhibited above 
remains intriguing and raises the hope that some of the structural aspects 
discussed in this contribution are generic features of quantum field theories.

 \vskip 4em


\small

\begin{thebibliography}{FjFRS2} \addtolength{\itemsep}{-4pt}
\addcontentsline{toc}{section}{References}

   \def\atmp  {Adv.\wb Theor.\wb Math.\wb Phys.}
   \def\bams  {Bull.\wb Amer.\wb Math.\wb Soc.}  
   \def\coma  {Con\-temp.\wb Math.}
   \def\comp  {Com\-mun.\wb Math.\wb Phys.}
   \def\cpma  {Com\-pos.\wb Math.} 
   \def\jgap  {J.\wb Geom.\wB and\wB Phys.}
   \def\jodg  {J.\wb Diff.\wb Geom.}
   \def\jopa  {J.\wb Phys.\ A}
   \def\nupb  {Nucl.\wb Phys.\ B}
   \def\phlb  {Phys.\wb Lett.\ B}
   \def\phrl  {Phys.\wb Rev.\wb Lett.}
   \def\pnas  {Proc.\wb Natl.\wb Acad.\wb Sci.\wb USA}
   \def\rvmp  {Rev.\wb Math.\wb Phys.}
   \def\taac  {Theo\-ry\wB and\wB Appl.\wb Cat.}
\newcommand\Bi[2]    {\bibitem[#1]{#2}}
\newcommand\BOOK[4]  {{\sl #1\/} ({#2}, {#3} {#4})}
\newcommand\inBO[7]  {{\sl #7\/}, in:\ {\sl #1}, {#2}\ ({#3}, {#4} {#5}), p.\ {#6}}
\newcommand\J[5]     {{\sl #5\/}, {#1} {#2} ({#3}) {#4}}
\newcommand\Jj[4]    {{\sl #4\/}, {#1} ({#2}) {#3}}
\newcommand\K[6]     {{\sl #6\/}, {#1} {#2} ({#3}) {#4}\,
                     \href{http://xxx.lanl.gov/abs/#5}{[{\tt #5}]}}
\newcommand\Prep[2]  {{\sl #2\/}, preprint {{\tt #1}}}
\def\wb{\,\linebreak[0]} 
\def\wB{\,\wb}

\Bi{Al}{alvaO5}{O.\ Alvarez, \J\comp{100}{1985}{279}
           {Topological quantization and cohomology}}
\Bi{ADW}{axdw}  {S.\ Axelrod, S.\ Della Pietra, and E.\ Witten, 
           \J\jodg{33}{1991}{787} {Geometric quantization of Chern-Simons theory}}
\Bi{BK}{BAki}   {B.\ Bakalov and A.A.\ Kirillov, \BOOK{Lectures on Tensor Categories
           and Modular Functors} {American Mathematical Society}{Providence}{2001}}
\Bi{FFFS}{fffs3} {G.\ Felder, J.\ Fr\"ohlich, J.\ Fuchs, and C.\ Schweigert,
           \K\cpma{131}{2002}{189} 
           {hep-th/9912239}
           {Correlation functions and boundary conditions in RCFT and three-dimensional topology}}
\Bi{FjFRS1}{fjfrs} {J.\ Fjelstad, J.\ Fuchs, I.\ Runkel, and C.\ Schweigert,
           \K\taac{16}{2006}{342} {hep-th/0503194} {TFT construction
           of RCFT corre\-la\-tors V: Proof of modular invariance and factorisation}} 
\Bi{FjFRS2}{fjfrs2}{J.\ Fjelstad, J.\ Fuchs, I.\ Runkel, and C.\ Schweigert,
           \K\atmp{12}{2008}{1281} 
           {hep-th/0612306}
           {Uniqueness of open/closed rational CFT with given algebra of open states}}
\Bi{FHT}{frht5} {D.S.\ Freed, M.J.\ Hopkins, and C.\ Teleman,
           \Prep{math.AT/0711.1906} {Loop groups and twisted K-theory I}}
\Bi{FrB}{FRbe2} {E.\ Frenkel and D.\ Ben-Zvi, \BOOK{Vertex Algebras and Algebraic Curves{\rm,
           second edition}} {American Mathematical Society}{Providence}{2004}}
\Bi{FRS1}{fuRs4} {J.\ Fuchs, I.\ Runkel, and C.\ Schweigert,
           \K\nupb{646}{2002}{353} {hep-th/0204148}
           {TFT construction of RCFT correlators I: Partition functions}}
\Bi{FRS2}{fuRs8} {J.\ Fuchs, I.\ Runkel, and C.\ Schweigert, \K\nupb{678}{2004}{511} 
           {hep-th/0306164} {TFT construction of RCFT correlators II: Unoriented world sheets}}
\Bi{FRS3}{fuRs10}{J.\ Fuchs, I.\ Runkel, and C.\ Schweigert,
           \K\nupb{715}{2005}{539} {hep-th/0412290} {TFT construction of 
           RCFT correlators IV: Structure constants and correlation functions}}
\Bi{FRS4}{fuRs11} {J.\ Fuchs, I.\ Runkel and C.\ Schweigert, \K\coma{431}{2007}{203} 
           {math.CT/0511590}
           {Ribbon categories and (unoriented) CFT: Frobenius algebras, automorphisms, reversions}}
\Bi{FSW}{fusW}  {J.\ Fuchs, C.\ Schweigert, and K.\ Waldorf, \K\jgap{58}{2008}{576}
           {hep-th/0703145} {Bi-branes: Target space geometry for world sheet topological defects}}
\Bi{Ga1}{gawe1}  {K.\ Gaw\c edzki, \inBO{Nonperturbative Quantum Field Theory} {G.\ 't Hooft,
           A.\ Jaffe, G.\ Mack, P.K.\ Mitter, and R.\ Stora, eds.} {Plenum}{New York}{1988} {101}
           {Topological actions in two-dimensional quantum field theories}}
\Bi{Ga2}{gawedzki4} {K.\ Gaw\c edzki, \K\comp{258}{2005}{23} {hep-th/0406072}
           {Abelian and non-Abelian branes in WZW models and gerbes}}
\Bi{GR}{gawedzki1}  {K.\ Gaw\c edzki and N.\ Reis,
           \K\rvmp{14}{2002}{1281} {hep-th/0205233} {WZW branes and gerbes}}
\Bi{GSW1}{gasuwa1}  {K.\ Gaw\c edzki, R.\ R.\ Suszek and K.\ Waldorf, \K\comp{284}{2008}{1}
           {hep-th/0701071} {WZW orientifolds and finite group cohomology}}
\Bi{GSW2}{gasuwa2}  {K.\ Gaw\c edzki, R.\ R.\ Suszek and K.\ Waldorf,
           \Prep{0809.5125} {Bundle gerbes for orientifold sigma models}}
\Bi{Hu}{huan21} {Y.-Z.\ Huang, \K\pnas{102}{2005}{5352} 
           {math.QA/0412261}
           {Vertex operator algebras, the Verlinde conjecture and modular tensor categories}}
\Bi{JW}{jewe}  {L.C.\ Jeffrey and J.\ Weitsman, \J\comp{150}{1992}{593} {Bohr-Sommerfeld orbits
           in the moduli space of flat connections and the Verlinde dimension formula}}
\Bi{KRT}{KArt} {C.\ Kassel, M.\ Rosso, and V.\ Turaev, \BOOK{Quantum Groups and Knot Invariants}
           {Soc.\ Math.\ de France}{Paris}{1997}}
\Bi{RS}{ruSu}  {I.\ Runkel and R.R.\ Suszek,
           \Prep{0808.1419} {Gerbe-holonomy for surfaces with defect networks}} 
\Bi{SaS}{saSc}  {G.\ Sarkissian and C.\ Schweigert, \Prep{0810.3159}
           {Some remarks on defects and T-duality}}
\Bi{SSW}{scsW}  {U.\ Schreiber, C.\ Schweigert, and K.\ Waldorf, \K\comp{274}{2007}{31} 
           {hep-th/0512283} {Unoriented WZW models and holonomy of bundle gerbes}}
\Bi{SFR}{scfr2} {C.\ Schweigert, J.\ Fuchs, and I.\ Runkel, \inBO{{\rm Proceedings of the}
           International Congress of Mathematicians 2006} {M.\ Sanz-Sol\'e, J.\ Soria, J.L.\ Varona,
           and J.\ Verdera, eds.} {European Mathematical Society}{Z\"urich}{2006} {443}
           {Categorification and correlation functions in conformal field theory}}
           \href{http://xxx.lanl.gov/abs/math.CT/0602079}{[{\tt math.CT/0602079}]}
\Bi{Tu}{TUra}  {V.G.\ Turaev, \BOOK{Quantum Invariants of Knots and $3$-Manifolds}
           {de Gruyter}{New York}{1994}}
\Bi{Wa1}{waldorf1} {K.\ Waldorf, \K\taac{17}{2007}{240}
           {math.CT/0702652} {More morphisms between bundle gerbes}}
\Bi{Wa2}{waldorf2} {K.\ Waldorf, \Prep{0804.4835}
           {Multiplicative bundle gerbes with connection}}
\end{thebibliography}
\end{document}